\input ./style/arxiv-general.cfg
\documentclass[aos,MSNbibl,seceqn,nameyear,dvips]{arximspdf}
\makeatletter
   \@ifpackageloaded{graphicx}{}{\usepackage{graphicx}}
\makeatother
\usepackage{mathbh}

\doi{10.1214/15-AOS1357}
\volume{43}
\issue{6}
\pubyear{2015}
\firstpage{2706}
\lastpage{2737}
\docsubty{FLA}

\makeatletter
\def\boldss{}

\newcommand{\rrvert}{\vert}

\newcommand{\rrVert}{\Vert}
\newcommand{\llvert}{\vert}
\newcommand{\llVert}{\Vert}
\newcommand{\eqref}[1]{(\ref{#1})}

\newcommand{\bS}{\mathbf{S}}
\newcommand{\bX}{\mathbf{X}}
\newcommand{\bone}{\mathbf{1}}

\newcommand{\tr}{\operatorname{tr}}
\newcommand{\diag}{\operatorname{diag}}
\newcommand{\offdiag}{\mathrm{off}}

\newcommand{\var}{\operatorname{var}}

\newtheorem{thmm}{Theorem}[section]
\newproclaim{defn}{Definition}[section]
\newtheorem{lem}{Lemma}[section]

\newtheorem{prop}{Proposition}[section]
\newproclaim{assum}{Assumption}[section]
\newproclaim{exm}{Example}[section]
\newproclaim{remark}{Remark}[section]

\newcommand{\KL}{\mathsf{KL}}

\makeatother

\begin{document}
\begin{frontmatter}

\title{Estimation of functionals of sparse covariance matrices}
\runtitle{Functionals of covariance matrices}

\begin{aug}
\author[A]{\fnms{Jianqing}~\snm{Fan}\thanksref{T1,m1}\ead[label=jqfan]{jqfan@princeton.edu}},
\author[B]{\fnms{Philippe}~\snm{Rigollet}\thanksref{T2,m2}\ead[label=rigollet]{rigollet@math.mit.edu}}
\and
\author[A]{\fnms{Weichen}~\snm{Wang}\corref{}\thanksref{T3,m1}\ead[label=weichenw]{weichenw@princeton.edu}}
\runauthor{J. Fan, P. Rigollet and W. Wang}
\affiliation{Princeton University\thanksmark{m1} and Massachusetts
Institute of Technology\thanksmark{m2}}

\thankstext{T1}{Supported in part by NSF Grants DMS-12-06464 and
DMS-14-06266 and NIH grant R01-GM072611-9.}
\thankstext{T2}{Supported in part by NSF Grants DMS-13-17308,
CAREER-DMS-1053987 and by the Howard~B. Wentz Jr. Junior Faculty award.}
\thankstext{T3}{Supported in part by NSF Grant DMS-12-06464.}

\address[A]{J. Fan\\
W. Wang\\
Department of Operations Research\\
\quad and Financial Engineering\\
Princeton University\\
Princeton, New Jersey 08544\\
USA\\
\printead{jqfan}\\
\phantom{E-mail:\ }\printead*{weichenw}}

\address[B]{P. Rigollet\\
Department of Mathematics\\
Massachusetts Institute of Technology\\
77 Massachusetts Avenue\\
Cambridge, Massachusetts 02139-4307\\
USA\\
\printead{rigollet}}

\end{aug}

%
\received{\smonth{2} \syear{2015}}
%
\revised{\smonth{6} \syear{2015}}


\begin{abstract}
High-dimensional statistical tests often ignore correlations to gain
simplicity and stability leading to null distributions that depend on
functionals of correlation matrices such as their Frobenius norm and
other $\ell_r$ norms. Motivated by the computation of critical values
of such tests, we investigate the difficulty of estimation the
functionals of sparse correlation matrices.
Specifically, we show that simple plug-in procedures based on
thresholded estimators of correlation matrices are sparsity-adaptive
and minimax optimal over a large class of correlation matrices. Akin to
previous results on functional estimation, the minimax rates exhibit an
elbow phenomenon.
Our results are further illustrated in simulated data as well as an
empirical study of data arising in financial econometrics.
\end{abstract}

%
\end{frontmatter}

\section{Introduction}
\label{SEC:intro}
Covariance matrices are at the core of many statistical procedures such
as principal component analysis or linear discriminant analysis.
Moreover, not only do they arise as natural quantities to capture
interactions between variables but, as we illustrate below, they often
characterize the asymptotic variance of commonly used estimators.
Following the original papers of \citeauthor{BicLev08a}
(\citeyear{BicLev08a,BicLev08b}), much work
has focused on the inference of high-dimensional covariance matrices
under sparsity [\citet
{CaiLiu11,CaiRenZho13,CaiYua12,CaiZhaZho10,CaiZho12,Kar08,LamFan09,RavWaiRas11}] and other structural assumptions related to sparse
principal component analysis [\citet{AmiWai09},
\citeauthor{BerRig13a} (\citeyear{BerRig13a,BerRig13b}),
\citet{BirJohNad13}, \citeauthor{CaiMaWu13} (\citeyear{CaiMaWu13,CaiMaWu14}),
\citet{JohLu09,LevVer12,RotLevZhu09,Ma13,OnaMorHal13,PauJoh12,FanFanLv08},
\citeauthor{FanLiaMin11} (\citeyear{FanLiaMin11,FanLiaMin13}),
\citet{SheSheMar11,VuLei12,ZouHasTib06}].
This area of research is very active and, as a result, this list of
references is illustrative rather than comprehensive. This line of work
can be split into two main themes: estimation and detection. The former
is the main focus of the present paper. However, while most of the
literature has focused on estimating the covariance matrix itself,
under various performance measures, we depart from this line of work by
focusing on functionals of the covariance matrix rather than the
covariance matrix itself.

Estimation of functionals of unknown signals such as regression
functions or densities is known to be different in nature from
estimation of the signal itself. This problem has received most
attention in nonparametric estimation, originally in the Gaussian white
noise model [\citet{IbrNemKha87,NemKha87,Fan91,EfrLow96}] [see also
\citet{Nem00} for a survey of results in the Gaussian white noise model]
and later extended to density estimation [\citet{HalMar87,BicRit88}] and
various other models such as regression [\citet{DonNus90}, \citeauthor{CaiLow05}
(\citeyear{CaiLow05,CaiLow06}),
\citet{Kle06}] and inverse problems [\citet{But07,ButMez11}]. Most of
these papers study the estimation of quadratic functionals and,
interestingly, exhibit an elbow in the rates of convergence: there
exists a critical regularity parameter below which the rate of
estimation is nonparametric and above which, it becomes parametric. As
we will see below the phenomenon also arises when regularity is
measured by sparsity.

Over the past decade, sparsity has become the prime measure of
regularity, both for its flexibility and generality. In particular,
smooth functions can be viewed as functions with a sparse expansion in
an appropriate basis. At a high level, sparsity assumes that many of
the unknown parameters are equal to zero or nearly so, so that the few
nonzero parameters can be consistently estimated using a small number
of observations relative to the apparent dimensionality of the problem.
Moreover, sparsity acts not only as a regularity parameter that
stabilizes statistical procedures but also as key feature for
interpretability. Indeed, it is often the case that setting many
parameters to zero simply corresponds to a simpler sub-model. The main
idea is to let data select the correct sub-model. This is the case in
particular for covariance matrix estimation where zeros in the matrix
correspond to uncorrelated variables. Yet, while the value of sparsity
for covariance matrix estimation has been well established, to the best
of our knowledge, this paper provides the first analysis for the
estimation of functionals of sparse covariance matrix. Indeed, the
actual performance of many estimators critically depends on such
functionals. Therefore, accurate functional estimation leads to a
better understanding the performance of many estimators and can
ultimately serve as a guide to selecting the best estimator.
Applications of our results are illustrated in Section~\ref{SEC:appli}.

Our work is not only motivated by real applications, but also by a
natural extension of the theoretical analysis carried out in the sparse
Gaussian sequence model [\citet{CaiLow05}]. In that paper, Cai and Low
assume that the unknown parameter $\theta$ belong to an $\ell_q$-ball,
where $q>0$ can be arbitrarily close to 0. Such balls are known to
emulate sparsity and actually correspond to a more accurate notion of
sparsity for signal $\theta$ that is encountered in applications [see,
e.g., \citet{FouRau13}]. They also show that a nonquadratic estimator
can be fully efficient to estimate quadratic functionals. 
We extend some of these results to covariance matrix estimation. Such
an extension is not trivial since, unlike the Gaussian sequence model,
covariance matrix lies at high-dimensional manifolds and its estimation
exhibits complicated dependencies in the structure of the noise.

We also compare our results for optimal rates of estimating matrix
functionals with that of estimating matrix itself. Many methods have
been proposed to estimate covariance matrix in different sense of
sparsity using different techniques including thresholding [\citet
{BicLev08a}], tapering [\citet{BicLev08b,CaiZhaZho10,CaiZho12}] and
penalized likelihood [\citet{LamFan09}] to name only a few. These methods
often lead to minimax optimal rates in various classes and under
several metrics [\citet{CaiZhaZho10,CaiZho12,RigTsy12a}]. However, the
optimal rates of estimating matrix functionals have not yet been
covered by much literature. Intuitively, it should have faster rates of
convergence on estimating a matrix functional than itself since it is
just a one-dimensional estimating problem and the estimating error
cancel with each other when we sum those elements together. We will see
this is indeed the case when we compare the minimax rates of estimating
matrix functionals with those of estimating matrices.


The rest of the paper is organized as follows. We begin in Section~\ref
{SEC:appli} by two motivating examples of high-dimensional hypothesis
testing problems: a two-sample testing problem of Gaussian means that
arises in genomics and validating the efficiency of markets based on
the Capital Asset Pricing Model (CAPM). Next, in Section~\ref
{SEC:quad}, we introduce an estimator of the quadratic functional of
interest that is based on the thresholding estimator introduced
in \citet
{BicLev08a}. We also prove its optimality in a minimax sense over a
large class of sparse covariance matrices. The study is further
extended to estimating other measures of sparsity of covariance matrix.
Finally, we study the numerical performance of our estimator in
Section~\ref{SEC:num} on simulated experiments as well as in the
framework of the two applications described in Section~\ref{SEC:appli}.
Due to space restrictions, the proofs for the upper bounds are
relegated to the \hyperref[app]{Appendix} in the supplementary
material [\citet{FanRigWan15}].

\textit{Notation}: Let $d$ be a positive integer. The space of $d
\times d$ positive semi-definite matrices is denoted by $\bS_d^+$. For
any two integers $c < d$, define $[c:d]=\{c, c+1, \ldots, d\}$ to be
the sequence of contiguous integers between $c$ and $d$, and we simply
write $[d]=\{1, \ldots, d\}$. $I_d$ denotes the identity matrix of
$\mathbb{R}^d$. Moreover, for any subset $S \subset[d]$, denote by
$\bone_S \in\{0,1\}^d$ the column vector with $j$th coordinate equal
to one iff $j \in S$. In particular, $\bone_{[d]}$ denotes the $d$
dimensional vector of all ones.

We denote by $\tr$ the trace operator on square matrices and by $\diag$
(resp., $\offdiag$) the linear operator that sets to 0 all the off
diagonal (resp., diagonal) elements of a square matrix. The Frobenius
norm of a real matrix $M$ is denoted by $\|M\|_\mathsf{F}$ and is
defined by
$\|M\|_\mathsf{F}=\sqrt{\tr(M^\top M)}$. Note that $\|M\|_\mathsf
{F}$ is a the
Hilbert--Schimdt norm associated with the inner product $\langle A,
B\rangle=\tr(A^\top B)$ defined on the space of real rectangular
matrices of the same size. Moreover, $|A|$ denotes the determinant of a
square matrix $A$. The variance of a random variable $X$ is denote by
$\var(X)$.

In the proofs, we often employ $C$ to denote a generic positive
constant that may change from line to line.

\section{Two motivating examples}
\label{SEC:appli}
In this section, we describe our main motivation for estimating
quadratic functionals of a high-dimensional covariance matrix in the
light of two applications to high-dimensional testing problems. The
first one is a high-dimensional two-sample hypothesis testing with
applications in gene-set testing. The second example is about testing
the validity of the \emph{capital asset pricing model (CAPM)} from
financial economics.

\subsection{Two-sample hypothesis testing in high-dimensions}
In various statistical applications, in particular in genomics, the
dimensionality of the problems is so large that statistical procedures
involving inverse covariance matrices are not viable due to its lack of
stability both from a statistical and numerical point of view. This
limitation can be well illustrated on a showcase example: two-sample
hypothesis testing [\citet{BaiSar96}] in high-dimensions.

Suppose that we observe two independent samples $X_1^{(1)}, \ldots,
X_{n_1}^{(1)} \in\mathbb{R}^p$ that are i.i.d. $\mathcal{N}(\mu
_1,\Sigma_1)$ and $X_1^{(2)}, \ldots, X_{n_2}^{(2)} \in\mathbb{R}^p$
that are i.i.d. $\mathcal{N}(\mu_2,\Sigma_2)$. Let $n = n_1 + n_2$. The
goal is to test $H_0: \mu_1 = \mu_2$ vs. $H_1: \mu_1 \ne\mu_2$.

Assume first that $\Sigma_1=\Sigma_2=\Sigma$. In this case, Hotelling's
test is commonly employed when $p$ is small.
Nevertheless, when $p$ is large,
\citet{BaiSar96} showed that the test based on Hotelling's $T^2$ has low
power and suggest a new statistics $M$ for the random matrix asymptotic
regime where $n, p \to\infty, \frac{n}{p} \to\gamma>0, \frac
{n_1}{n_1+n_2} \to\kappa\in(0,1)$. The statistics, implementing the
naive Bayes rule, is defined as
\[
M = \bigl(\bar X^{(1)} - \bar X^{(2)}\bigr)^{\top}
\bigl(\bar X^{(1)} - \bar X^{(2)}\bigr) - \frac{n}{n_1n_2} \tr(
\hat\Sigma) ,
\]
and is proved to be asymptotically normal under the null hypothesis with
%
\[
\var(M) = 2\frac{n(n-1)}{(n_1n_2)^2}\|\Sigma\|_\mathsf{F}^2
\bigl(1+o(1)\bigr).
\]
Clearly, the asymptotic variance of $M$ depends on the unknown
covariance matrix $\Sigma$ through its quadratic functional, and in
order to compute the critical value of the test, Bai and Saranadasa
suggest to estimate $\|\Sigma\|_\mathsf{F}^2$ by the quantity
\[
B^2 = \frac{n^2}{(n+2)(n-1)} \biggl[\|\hat\Sigma\|_\mathsf{F}^2-
\frac
{1}{n}\bigl(\tr (\hat\Sigma)\bigr)^2 \biggr].
\]
They show that $B^2$ is a ratio-consistent estimator of $\|\Sigma \|
_\mathsf{F}^2$ in the sense that $B^2=(1+o_P(1))\|\Sigma\|_\mathsf
{F}^2$. Clearly, this
solution does not leverage any sparsity assumption and may suffer from
power deficiency if the matrix $\Sigma$ is indeed sparse. Rather, if
the covariance matrix $\Sigma$ is believed to be sparse, one may prefer
to use a thresholded estimator for $\Sigma$ as in \citet{BicLev08a}
rather than the empirical covariance matrix $\hat\Sigma$. In this case,
we estimate $\|\Sigma\|_\mathsf{F}^2$ by $\widehat{\|\Sigma\|
_\mathsf{F}^2} =
\sum_{i,j = 1}^p {\hat\sigma_{ij}^2} \mathbh{1} \{ {\llvert  {{{\hat
\sigma}_{ij}}} \rrvert  > \tau}  \}$, where $\{\hat\sigma
_{ij},
i, j \in[p]\}$ could be any consistent estimator of $\sigma_{ij}$ and
$\tau>0$ is a threshold parameter.

More recently, \citet{CheQin10} took into account the case $\Sigma_1
\neq\Sigma_2$ and proposed a test statistic based on an unbiased
estimate of each of the three quantities in $\|\mu_1 -\mu_2\|^2 =\|
\mu
_1\|^2 + \|\mu_2\|^2 - 2 \mu_1^{\top}\mu_2$. 
In this case, the quantities $\|\Sigma_i\|_\mathsf{F}^2$, $i=1,2$ and
$\langle
\Sigma_1, \Sigma_2\rangle$ appear in the asymptotic variance. The
detailed formulation and assumptions of this statistic, as well as
discussions about other testing methods such as \citet{SriDu08}, are
provided in the supplementary material [\citet{FanRigWan15}] for completeness. If
$\Sigma_1$ and $\Sigma_2$ are indeed sparse, akin to the above
reasoning, we can also estimate $\|\Sigma_i\|_\mathsf{F}^2$, $i=1,2$ and
$\langle\Sigma_1, \Sigma_2\rangle$ using thresholding to leverage
sparsity assumption. It is not hard to derive a theory for estimating
quadratic functionals involving two covariance matrices but the details
of this procedure are beyond the scope of the present paper.
\subsection{Testing high-dimensional CAPM model}
The capital asset pricing model (CAPM) is a simple financial model that
postulates how individual asset returns are related to the market
risks. Specifically, the individual excessive return $Y^{(i)}_t$ of asset
$i \in[N]$ over the risk-free rate at time $t \in[T]$ can be
expressed as an affine function of a vector of $K$ risk factors $f_t
\in\mathbb{R}^K$:
%
\begin{equation}
\label{EQ:capm} Y^{(i)}_t = \alpha_i +
\beta_i^\top f_t + \varepsilon^{(i)}_{t}
,
\end{equation}
where we assume for any $t \in[T]$, $f_t \in\mathbb{R}^K$ are
observed. The case $K=1$ with $f_t$ being the excessive return of the
market portfolio corresponds to the CAPM [\citet{Sha64,Lin65,Mos66}].
It is nowadays more common to employ the Fama--French three-factor
model [see \citet{FamFre93} for a definition] for the US equity market,
corresponding to $K=3$.

For simplicity, let us rewrite the model \eqref{EQ:capm} in the
vectorial form
\[
Y_t = \alpha+ Bf_t + \varepsilon_t , \qquad t
\in[T].
\]
%
The multi-factor pricing model 
postulates $\alpha= 0$. Namely, all returns are fully compensated by
their risks: no extra returns are possible and the market is efficient.
This leads us to naturally consider the hypothesis testing problem
$H_0: \alpha= 0$ vs. $H_1: \alpha\ne0$.




Let $\hat{\alpha}$ and $\hat{B}$ be the least-squares estimate and
$\hat{\varepsilon}_t = Y_t - \hat{\alpha} - \hat{B} f_t$ be a
residual vector.
Then an unbiased estimator of $\Sigma= \var(\varepsilon_t)$ is
\[
\tilde\Sigma=\frac{1}{T-K-1}\sum_{t=1}^T
\hat\varepsilon_t\hat \varepsilon _t^\top.
\]
Let $\hat D=\diag(\tilde\Sigma)$ and $M_F=I_T-F(F^\top F)^{-1}F^\top$
where $F = (f_1,\ldots, f_T)^{\top}$. Define $W_d = (\bone
_{[T]}^\top
M_F \bone_{[T]}){ \hat\alpha}^\top{\hat D}^{-1} {\hat\alpha}$ the
Wald-type of test statistics with correlation ignored, whose normalized
version is given by
%
\begin{equation}
\label{eq2.3} J_{\alpha} = \frac{W_d - \mathbb{E}(W_d)}{\sqrt{\var(W_d)}}.
\end{equation}



Under some conditions, it was shown by \citet{PesYam12} that, under
$H_0$, $J_\alpha\to\mathcal{N}(0,1)$ as $N \to\infty$. Moreover, if
$\varepsilon^{(i)}_t$'s are i.i.d. Gaussian, it holds that $\mathbb
{E}(W_d) =
\nu N/(\nu-2)$ and
\[
\var(W_d) = \frac{2N(\nu-1)}{\nu-4} \biggl(\frac{\nu}{\nu-2}
\biggr)^2 \bigl[1+(N-1) {\bar{\rho}^2}+O\bigl(
\nu^{-1/2}\bigr)\bigr] ,
\]
where $\nu=T-K-1$ is the degrees of freedom and
\[
{\bar{\rho}^2} = \frac{2}{N(N-1)} \sum
_{i=2}^N \sum_{j=1}^{i-1}
\rho _{ij}^2 ,
\]
where $\rho=D^{-1/2}\Sigma D^{-1/2}$ with $D=\diag(\Sigma)$ is the
correlation matrix of the stationary process $(\varepsilon_t)_{t\in[T]}$.
The authors go on to propose an estimator of the quadratic functional
${\bar{\rho}^2}$ by replacing the correlation coefficients $\rho_{ij}$
in the above expression by $\hat\rho_{i,j}\mathbh{1}(|\hat\rho
_{ij}|>\tau)$ where $(\hat\rho_{ij})_{i,j \in[N]}=\hat
D^{-1/2}\tilde
\Sigma\hat D^{-1/2}$ and $\tau>0$ is a threshold parameter. However,
they did not provide any analysis of this method, nor any guidance to
chose $\tau$. 

\section{Optimal estimation of quadratic functionals}
\label{SEC:quad}

In the previous section, we have described rather general questions
involving the estimation of quad\-ratic functions of covariance or
correlation matrices. We begin by observing that consistent estimation
of $\|\Sigma\|_\mathsf{F}^2$ is impossible unless $p=o(n)$. This
precludes in
particular the high-dimensional framework that motivates our study.

Our goal is to estimate the Frobenius norm $\|\Sigma\|_\mathsf{F}^2$
of a
sparse $p \times p$ covariance matrix $\Sigma$ using $n$ i.i.d.
observations $X_1, \ldots, X_n \sim\mathcal{N}(0, \Sigma)$. Observe
that $\|\Sigma\|_\mathsf{F}^2$ can be decomposed as $\|\Sigma \|
_\mathsf{F}^2=Q(\Sigma
)+D(\Sigma)$ where $Q(\Sigma)= \sum_{i \ne j} {\sigma_{ij}^2}$
corresponds to the off-diagonal elements and $D(\Sigma)=\sum_j
{\sigma
_{jj}^2}$ corresponds to the diagonal elements. The following theorem,
implies that even if $\Sigma=\diag(\Sigma)$ is diagonal, the quadratic
functional $\|\Sigma\|_\mathsf{F}^2$ cannot be estimated consistently in
absolute error if $p\ge n$. Note that the situation is quite different
when it comes to \emph{relative error}. Indeed, the estimator of \citet
{BaiSar96} is consistent in relative error with no sparsity assumption
even in the high-dimensional regime. Study of the relative error in the
presence of sparsity is an interesting question that deserves further
developments.
This makes sense intuitively as the diagonal of $\Sigma$ consists of
$p$ unknown parameters while we have only $n$ observations.

\begin{prop}
\label{PROP:diag}
Fix $n, p \ge1$ and let
\[
\mathcal{D}_p=\bigl\{\Sigma\in\bS_p^+: \Sigma=\diag(
\Sigma) , \Sigma_{ii} \le1 \bigr\}
\]
be the class of diagonal covariance matrices with diagonal elements
bounded by~1. Then there exists a universal constant $C>0$ such that
\[
\inf_{\hat D} \sup_{\Sigma\in\mathcal{D}_p} \mathbb{E} \bigl[\hat
D- D(\Sigma) \bigr]^2 \ge C {\frac{p}{n}}.
\]
In particular, it implies that
\[
\inf_{\hat F} \sup_{\Sigma\in\mathcal{D}_p} \mathbb{E} \bigl[\hat
F - \|\Sigma\|_\mathsf{F}^2 \bigr]^2 \ge C{
\frac{p}{n}},
\]
where the infima are taken with over all real valued measurable
functions of the observations.
\end{prop}

\begin{pf}
Our lower bounds rely on standard arguments from minimax theory. We
refer to Chapter~2 of \citet{Tsy09} for more details. In the sequel, let
$\KL(P, \bar P)$ denote the Kullback--Leibler divergence between two
distributions $P$ and $\bar P$, where $P \ll\bar P$. It is defined by
\[
\KL(P,\bar P)=\int\log \biggl(\frac{\mathrm{d} P}{\mathrm{d} \bar
P} \biggr)\,\mathrm{d} P.
\]

We are going to employ a simple two-point lower bound. Fix $\varepsilon
\in(0,1/2)$ and let $P_p^n$ (resp., $\bar P_p^n$) denote the
distribution of a sample $X_1, \ldots, X_n$ where $X_1 \sim\mathcal
{N}(0, I_p)$ [resp., $X_1 \sim\mathcal{N}(0, (1-\varepsilon) I_p)$].
Next, observe that $I_p, (1-\varepsilon)I_p \subset\mathcal{D}_p$ so that
%
\begin{equation}
\label{EQ:pf:diag1} \sup_{\Sigma\in\mathcal{D}_p} \mathbb{E}
\bigl|\hat D - D(\Sigma ) \bigr| \ge
\max_{\Sigma\in\{I_p, (1-\varepsilon)I_p\}} \mathbb{E}
 \bigl|\hat D - D(\Sigma)\bigr |.
\end{equation}
Moreover, $ |D(I_p)-D ((1-\varepsilon)I_p ) |=p(2\varepsilon
-\varepsilon^2)> p\varepsilon$. Then it follows from the Markov inequality that
%
\begin{eqnarray}
\label{EQ:pf:diag2}
\frac{1}{p\varepsilon}\max_{\Sigma\in\{I_p, (1-\varepsilon)I_p\}} \mathbb {E} \bigl|\hat D
- D(\Sigma)\bigr | & \ge& \max_{\Sigma\in\{I_p,
(1-\varepsilon)I_p\}}\mathbb{P} \bigl[ \bigl|\hat D - D(
\Sigma) \bigr|>p \varepsilon \bigr]
\nonumber
\\[-8pt]
\\[-8pt]
\nonumber
& \ge& \frac{1}{4}\exp\bigl[-\KL\bigl(P_p^n, \bar
P_p^n\bigr)\bigr],
\end{eqnarray}
where the last inequality follows from Theorem~2.2(iii) of
\citet{Tsy09}.

Completion of the proof requires an upper bound on $\KL(P_p^n, \bar
P_p^n)$. To that end, note that it follows from the chain rule and
simple algebra that
\[
\KL\bigl(P_p^n, \bar P_p^n
\bigr)=np\KL\bigl(P_1^1, \bar P_1^1
\bigr)=\frac{np}{2} \biggl[\log (1-\varepsilon)+\frac{\varepsilon}{1-\varepsilon} \biggr]\le
\frac
{np}{2}\frac
{\varepsilon^2}{1-\varepsilon}\le np\varepsilon^2.
\]
Taking now $\varepsilon=1/(2\sqrt{np})\le1/2$ yields $\KL(P_p^n,
Q_p^n)\le1/4$. Together with \eqref{EQ:pf:diag1} and \eqref
{EQ:pf:diag2}, it yields
\[
\inf_{\hat D}\sup_{\Sigma\in\mathcal{D}_p} \mathbb{E}\bigl |\hat D - D(
\Sigma)\bigr | \ge\frac{1}{8e^{1/4}}\sqrt{\frac{p}{n}}.
\]
To complete the proof, we square the above inequality and employ
Jensen's inequality.
\end{pf}


To overcome the above limitation, we consider the following class of
sparse covariance matrices (indeed correlation matrices). For any $q
\in[0,2), R>0$ let $\mathcal{F}_q(R)$ denote the set of $p \times p$
covariance matrices defined by
%
\begin{equation}
\label{EQ:defFq} \mathcal{F}_q(R)= \biggl\{ \Sigma\in
\bS_p^+: \sum_{i \ne j} |\sigma
_{ij}|^q \le R , \diag(\Sigma)= I_p \biggr\}.
\end{equation}
Note that for this class of functions, we assume that the variance
along each coordinate is normalized to 1. This normalization is
frequently obtained by sample estimates, as shown in the previous
section. This simplified assumption is motivated also by
Proposition~\ref{PROP:diag} above which implies that $\|\Sigma\|
_\mathsf{F}^2$
for general covariance matrix cannot be estimated accurately in
absolute error in the large $p$ small $n$ regime since sparsity
assumptions on the diagonal elements are implausible. Note that the
condition $\diag(\Sigma)=I_p$ implies that diagonal elements
$D(\Sigma
)$ of matrices in $\mathcal{F}_q(R)$ can be estimated without error so
that we could possibly achieve consistency even if the case of large
$p$ small $n$.

Matrices in $\mathcal{F}_q(R)$ have many small coefficients for small
values of $q$ and~$R$. In particular, when $q=0$, there are no more
than $R$ entries of nonvanishing correlations. Following a major trend
in the estimation of sparse covariance matrices [\citeauthor{BicLev08a} (\citeyear{BicLev08a,BicLev08b}),
\citet{CaiLiu11,CaiYua12,CaiZhaZho10,CaiZho12,Kar08,LamFan09}], we employ a thresholding estimator of the covariance matrix
as a running horse to estimate the quadratic functionals. From the $n$
i.i.d. observations $X_1, \ldots, X_n \sim\mathcal{N}(0, \Sigma)$, we
form the empirical covariance matrix ${\hat\Sigma}$ that is defined by
%
\begin{equation}
\label{EQ:defempcov} {\hat\Sigma}=\frac{1}{n} \sum
_{k=1}^n X_k X_k^{\top}
\end{equation}
with elements $\hat\Sigma=\{\hat\sigma_{ij}\}_{ij}$ and for any
threshold $\tau>0$, let $\tilde\Sigma_\tau=\{\tilde\sigma_{ij}\}
_{ij}$ denote the thresholding estimator of $\Sigma$ defined by
$\tilde
\sigma_{ij}=\hat\sigma_{ij}\mathbh{1}\{|\hat\sigma_{ij}|>\tau\}$ if
$i\neq j$ and $\tilde\sigma_{ii}=\hat\sigma_{ii}$.

Next, we employ a simple plug-in estimator for $Q(\Sigma)$:
%
\begin{equation}
\label{EstQ} \widehat{Q(\Sigma)}=Q(\tilde\Sigma_\tau) = \sum
_{i \ne j} {\hat \sigma_{ij}^2}
\mathbh{1} \bigl\{ {\llvert {{{\hat\sigma}_{ij}}} \rrvert > \tau} \bigr
\} .
\end{equation}
Note that no value of the diagonal elements is used to estimate
$Q(\Sigma)$.

In the rest of this section, we establish that $\widehat{Q(\Sigma)}$ is
minimax adaptive over the scale $\{\mathcal{F}_q(R), q\in[0,2), R>0\}
$. Interestingly, we will see that the minimax rate presents an elbow
as often in quadratic functional estimation.

\begin{thmm}
\label{TH:UBquad}
Assume that $\gamma\log(p)<n$ for some constant $\gamma>8$ and fix
$C_0 \ge4$. Consider the threshold
\[
\tau= 2C_0\sqrt{\frac{\gamma\log p}{n}} ,
\]
and assume that $\tau\le1$.
Then, for any $q \in[0,2), R>0$, the plug-in estimator $Q(\tilde
\Sigma
_\tau)$ satisfies
\[
\mathbb{E} \bigl[ \bigl(Q(\tilde\Sigma_\tau) - Q(\Sigma)
\bigr)^2 \bigr] \le C_1\psi_{n,p}(q,R) +
C_2p^{4-\gamma/2} ,
\]
where
\[
\psi_{n,p}(q,R)= \frac{R^{2}}{n} \vee R^2 \biggl(
\frac{\log p}{n} \biggr)^{2-q} ,
\]
and $C_1$, $C_2$ are positive constants depending on $\gamma, C_0, q$.
\end{thmm}

The proof is postponed to the supplementary material.

Note that the rates $\psi_{n,p}(q,R)$ present an elbow at $q = 1 -
\log
\log p/\log n$ as usually the case in functional estimation.
We now argue that the rates $\psi_{n,p}(q,R)$ are optimal in a minimax
sense for a wide range of settings. In particular, the elbow effect
arising from the maximum in the definition of $\psi$ is not an
artifact. In the following theorem, we emphasize the dependence on
$\Sigma$ by using the notation $\mathbb{E}_\Sigma$ for the expectation
with respect to the distribution of the sample $X_1, \ldots, X_n$,
where $X_i \sim\mathcal{N}(0,\Sigma)$.

\begin{thmm}
\label{TH:LBquad}
Fix $q \in[0,2), R>0$ and assume $2\log p < n$ and $R^2 < (p-1)n^{-q}
/2$. Then there exists a positive constant $C_3>0$ such that
\[
\mathop{\inf} _{{{\hat Q}}} \mathop{\sup} _{\Sigma\in{\mathcal
F_q(R)}}
\mathbb{E}_{\Sigma} \bigl[ \bigl({\hat Q} - Q(\Sigma) \bigr)^2
\bigr]\ge C_3 \phi_{n,p}(q,R) ,
\]
where $\phi_{n,p}(q,R)$ is defined by
%
\begin{equation}
\label{EQ:phiquad} \phi_{n,p}(q,R)= \frac{R^2}{n} \boldss{\vee}
\biggl\{ R^2 \biggl(\frac
{\log ({(p-1)}/{(R^2n^{q})}+1 )}{2n} \biggr)^{2-q}\wedge
R^{4/q} \wedge1 \biggr\}
\end{equation}
and the infimum is taken over all measurable functions $\hat Q$ of the
sample $X_1, \ldots, X_n$.
\end{thmm}

Before proceeding to the proof, a few remarks are in order.
\begin{longlist}[1.]
\item[1.] The additional term of order $p^{4-\gamma/2}$ in Theorem~\ref
{TH:UBquad} can be made negligible by taking $\gamma$ large enough. To
show this tradeoff explicitly, we decided keep this term.
\item[2.] When $ 1\le R^2 < p^\alpha n^{-q}$ for some constant $\alpha<
1$, a slightly stronger requirement than Theorem~\ref{TH:LBquad},
the lower bound there can be written as
%
\begin{equation}
\label{eq:opt} \phi_{n,p}(q,R)= \frac{R^2}{n} \boldss{\vee}
\biggl\{ R^2 \biggl(\frac
{\log p}{n} \biggr)^{2-q} \wedge1
\biggr\}.
\end{equation}

Observe that the above lower bound matches the upper bound presented in
Theorem~\ref{TH:UBquad} when $R^{{2}/{(2-q)}} \log p \le n$.
Arguably, this is the most interesting range as it characterizes rates
of convergence (to zero) rather than rates of divergence, that may be
of different nature [see, e.g., \citet{Ver12}]. In other words, the
rates given in \eqref{eq:opt} are minimax adaptive with respect to $n$,
$R$, $p$ and $q$. In our formulation, we allow $R=R_{n,p}$ to depend
on other parameters of the problem. We choose here to keep the notation light.

%
\item[3.] The reason we choose correlation matrix class to present the
elbow effect is just for simplicity. Actually, we can replace the
constraint $\diag(\Sigma) = I_p$ in the definition of $\mathcal F_q(R)$
by boundedness of diagonal elements of $\Sigma$. Then for estimating
off-diagonal elements $Q(\Sigma)$, following exactly the same
derivation, the same elbow phenomenon has been noticed. Meanwhile, the
optimal rate for estimating diagonal elements $D(\Sigma)$ is again of
the order $p/n$. This optimal rate can be attained by the estimator
%
\begin{equation}
\label{EstD} \widehat{D(\Sigma)} = \frac{1}{n(n-1)} \sum
_{i=1}^p \sum_{k \ne j}
X_{k,i}^2 X_{j,i}^2.
\end{equation}
We omitted the proof here. Thus, if we do not have prior information
about diagonal elements, we could still estimate optimally the
quadratic functional of a covariance matrix by applying the
thresholding method \eqref{EstQ} for off-diagonal elements, together
with \eqref{EstD} for diagonal elements.
%
\item[4.] The rate $\phi_{n,p}(q,R)$ presents the same elbow phenomenon at
$q=1$ observed in the estimation of functionals, starting independently
with work of \citet{BicRit88} and \citet{Fan91}. Closer to the present
setup is the work of \citet{CaiLow05} who study the estimation of
functionals of ``sparse'' sequences in the infinite Gaussian sequence
model. There, a parameter controls the speed of decay of the unknown
coefficients. Note that while smaller values $q$ lead to sparser
matrices $\Sigma$, no estimator can benefit further from sparsity below
$q=1$ [the estimator has a rate of convergence $O(R^2/n)$ for any $q <
1$], unlike in the case of estimation of $\Sigma$. Again, this is
inherent to estimating functionals.
\item[5.] The condition $R^2 < (p-1)n^{-q}/2$ corresponds to the
high-dimensional regime and allows us to keep clean terms in the
logarithm. Similar assumptions are made in related literature [see,
e.g., \citet{CaiZho12}].
\item[6.] The optimal rates obtained here cannot be implied by existing
ones for estimating sparse covariance matrices. In particular, the
latter do not admit an elbow phenomenon. Specifically, \citet{RigTsy12a}
showed the optimal rate for estimating $\Sigma$ for $\Sigma\in
{\mathcal F_q}(R)$ under the Frobenius norm
is $\sqrt{R} {({\log p}/n)^{1/2 - q/4}}$ for $0 \le q < 2$.
Using this,
it is not hard to derive with high probability,
\[
\bigl|Q(\hat\Sigma) - Q(\Sigma)\bigr| \le C_1 R \biggl(\frac{\log p}{n}
\biggr)^{1/2 -
q/4} + C_2 R \biggl( \frac{\log p}{n}
\biggr)^{1-q/2} ,
\]
since $\| Q(\Sigma)\|_F = O(\sqrt{R})$ if nonvanishing correlations
are bounded away from zero. On one hand, when $q < 2$ the first term
always dominates so that we do not observe the elbow effect. In
addition, the rate so obtained is not optimal.
\end{longlist}
We now turn to the proof of Theorem~\ref{TH:LBquad}

\begin{pf*}{Proof of Theorem \ref{TH:LBquad}}
To prove minimax lower bounds, we employ a standard technique that
consists of reducing the estimation problem to a testing problem. We
split this proof into two parts and begin by proving
\[
\mathop{\inf} _{{{\hat Q}}} \mathop{\sup} _{\Sigma\in{\mathcal
F_q(R)}}
\mathbb{E}_{\Sigma} \bigl[{\hat Q} - Q(\Sigma) \bigr]^2\ge C
\frac
{R^{2}}{n} ,
\]
for some positive constant $C>0$. To that end, for any $A \in\bS^+_p$,
let $\mathbb{P}_A$ denote the distribution of $X \sim\mathcal
{N}(0,A)$. It is not hard to show if $|A|>0$ and $|B|>0$, $A,B \in\bS
^+_p$, then the Kullback--Leibler divergence between $\mathbb{P}_A$ and
$\mathbb{P}_B$ is given by
%
\begin{equation}
\label{fan1} \KL(\mathbb{P}_A, \mathbb{P}_B)=
\frac{1}{2} \biggl[\log \biggl(\frac
{|B|}{|A|} \biggr)+\tr
\bigl(B^{-1}A\bigr)-p \biggr].
\end{equation}
Next, take $A$ and $B$ to be of the form
\[
A^{(k)}=\pmatrix{ \bone
\bone^\top& a \bone\bone^\top& 0
\vspace*{2pt}\cr
a \bone \bone^\top& \bone\bone^\top& 0
\vspace*{2pt}\cr
0 & 0 & I_{p-k} } ,\qquad
B^{(k)}=\pmatrix{ \bone
\bone^\top& b \bone\bone^\top& 0
\vspace*{2pt}\cr
b \bone \bone^\top& \bone\bone^\top& 0
\vspace*{2pt}\cr
0 & 0 & I_{p-k} },
\]
where $a,b \in(0,1/2)$, $0$ is a generic symbol to indicate that the
missing space is filled with zeros, and $\bone$ denotes a vector of
ones of length $k/2$. Note that if we have random variables $(X, Y,
Z_1, \ldots, Z_{p-2})$ chosen from distribution $\mathcal{N}(0,
A^{(2)})$ meaning that $Z_k$'s are independent with $X, Y$ but the
correlation between $X$ and $Y$ is $a$, then random vector $(X, \ldots,
X, Y, \ldots, Y, Z_1, \ldots, Z_{p-k})$ with $k/2$ $X$'s and $Y$'s in it
follows $\mathcal{N}(0, A^{(k)})$. It is obvious that these two
matrices are degenerate and comes from perfectly correlated random
variables. Since perfectly correlated random variables do not add new
information, for such matrices, an application of \eqref{fan1} yields
\[
\KL(\mathbb{P}_{A^{(k)}}, \mathbb{P}_{B^{(k)}})=\KL(\mathbb
{P}_{A^{(2)}}, \mathbb{P}_{B^{(2)}}) = \frac{1-ab}{1-b^2}-
\frac
{1}{2}\log \biggl(\frac{1-a^2}{1-b^2} \biggr) -1.
\]
Next, using the convexity inequality $\log(1+x) \ge x-x^2/2$ for all
$x>0$, we get that
\[
\KL(\mathbb{P}_{A^{(k)}}, \mathbb{P}_{B^{(k)}})\le\frac
{(a-b)^2}{2(1-b^2)}
\biggl[1+\frac{(a+b)^2}{2(1-b^2)}\biggr]\le2(a-b)^2 ,
\]
using the fact that $a,b \in(0,1/2)$.
Take now if $R > 4$
\[
a=\frac{1}{4} ,\qquad b=a+\frac{1}{4 \sqrt{n}} ,\qquad k = \sqrt{R}
\]
so that we indeed have $a, b \in(0,1/2)$ and also $A^{(k)}, B^{(k)}
\in\mathcal{F}_q(R)$ obviously. If $R < 4$, take $k=2, a = \sqrt{R}/8,
b = a+\sqrt{R/64 n}$ instead. Moreover, this choice leads to $n\KL
(\mathbb{P}_A, \mathbb{P}_B)\le1/5$. Using standard techniques to
reduce estimation problems to testing problems [see, e.g., Theorem~2.5
of \citet{Tsy09}], we find that
\[
\inf_{\hat Q} \max_{\Sigma\in\{A,B\}}
\mathbb{E}_\Sigma \bigl[ \bigl(\hat Q-Q(\Sigma) \bigr)^2
\bigr]\ge C \bigl(Q(A)-Q(B) \bigr)^2.
\]
For the above choice of $A$ and $B$, we have
\[
\bigl(Q\bigl(A^{(k)}\bigr)-Q\bigl(B^{(k)}\bigr)
\bigr)^2 = \frac{k^4}{4} \bigl(a^2 - b^2
\bigr)^2 
\ge C \frac{R^2}{n}.
\]
Since $A^{(k)}, B^{(k)} \in\mathcal{F}_q(R)$, the above two displays
imply that
\[
\inf_{\hat Q} \max_{\Sigma\in\mathcal{F}_q(R)}
\mathbb{E}_\Sigma \bigl[ \bigl(\hat Q-Q(\Sigma) \bigr)^2
\bigr]\ge C \frac{R^2}{n} ,
\]
which completes the proof of the first part of the lower bound.

For the second part of the lower bound, we reduce our problem to a
testing problem of the same flavor as \citet{AriBubLug14,BerRig13b}.
Note, however, that our construction is different because the
covariance matrices considered in these papers do not yield large
enough lower bounds. We use the following construction.

Fix an integer $k \in[p-1]$ and let $\mathcal{S}=\{S \subset[p-1]:
|S|=k\}$ denote the set of subsets of $[p-1]$ that have cardinality
$k$. Fix $a\in(0, 1)$ to be chosen later and for any $S\in\mathcal
{S}$, recall that $\bone_S$ is the column vector in $\{0,1\}^{p-1}$
with support given by~$S$. For each $S \in\mathcal{S}$, we define the
following $p\times p$ covariance matrix:
%
\begin{equation}
\label{EQ:defSigmaS} \Sigma_S=\pmatrix{1
& a\bone_S^\top
\vspace*{2pt}\cr
a\bone_S & I_{p-1} }.
\end{equation}
Let $\mathbb{P}_0$ denote the distribution of $X \sim\mathcal{N}_p(0,
I_p)$ and $\mathbb{P}_S$ denote the distribution of $X \sim\mathcal
{N}_p(0, \Sigma_S)$. Let $\mathbb{P}_0^n$ (resp., $\mathbb{P}_S^n$)
denote the distribution of $\bX=(X_1, \ldots, X_n)$ of a collection $n$
i.i.d. random variables drawn from $\mathbb{P}_0$ (resp., $\mathbb
{P}_S$). Moreover, let $\bar{\mathbb{P}}^n$ denote the distribution of
$\bX$ where the $X_i$'s are drawn as follows: first draw $S$ uniformly
at random from $\mathcal{S}$ and then, conditionally on $S$, draw $X_1,
\ldots, X_n$ independently from $\mathbb{P}_S$. Note that $\bar
{\mathbb
{P}}^n$ is the mixture of $n$ independent samples rather the
distribution of $n$ independent random vectors drawn from a mixture
distribution.
Consider the following testing problem:
\[
H_0:\qquad \bX\sim\mathbb{P}_0^n \quad\mbox{vs.}\quad
H_1:\qquad \bX\sim\bar{\mathbb{P}}^n.
\]
Using Theorem~2.2, part (iii) of \citet{Tsy09}, we get that for any test
$\psi=\psi(\bX)$, we have
\[
\mathbb{P}_0^n(\psi=0)\vee\max_{S\in\mathcal{S}}
\mathbb{P}_S^n (\psi =1)\ge\mathbb{P}_0^n(
\psi=0)\vee\bar{\mathbb{P}}^n (\psi=1)\ge \frac{1}4\exp
\bigl(-\chi^2\bigl(\bar{\mathbb{P}}^n,
\mathbb{P}_0\bigr) \bigr) ,
\]
where we recall that the $\chi^2$-divergence between two probability
distributions $P$ and $Q$ is defined by
\[
\chi^2(P, Q)=\cases{ %
\displaystyle\int \biggl(
\frac{\mathrm{d} P}{\mathrm{d} Q}-1 \biggr)^2 \,\mathrm{d} Q , &\quad
 $\mbox{if } P\ll Q,$
\vspace*{2pt}\cr
\infty, &\quad $\mbox{otherwise.}$}
\]
Lemma~\ref{LEM:chi2div} implies that for suitable choices of the
parameters $a$ and $k$, we have $\chi^2(\bar{\mathbb{P}}^n, \mathbb
{P}_0)\le2$ so that the test errors are bounded below by a constant
$C=e^{-2}/4$. Since $Q(\Sigma_S)=2ka^2$ for any $S \in\mathcal{S}$, it
follows from a standard reduction from hypothesis testing to estimation
[see, e.g., Theorem~2.5 of \citet{Tsy09}] that the above result implies
the following lower bound:
%
\begin{equation}
\label{EQ:lbtheta} \inf_{\hat Q} \max_{\Sigma\in\mathcal{H}}
\mathbb{E}_\Sigma \bigl[ \bigl(\hat Q-Q(\Sigma) \bigr)^2
\bigr]\ge Ck^2a^4 ,
\end{equation}
for some positive constant $C$, where the infimum is taken over all
estimators $\hat Q$ of $Q(\Sigma)$ based on $n$ observations and
$\mathcal{H}$ is the class of covariance matrices defined by
\[
\mathcal{H}=\{I_p\}\cup\{\Sigma_S: S \in\mathcal{S}\}.
\]
To complete the proof, observe that the values of $a$ and $k$
prescribed in Lemma~\ref{LEM:chi2div} imply that $\mathcal{H} \subset
\mathcal{F}_q(R)$ and give the desired lower bound. Note first that,
for any choice of $a$ and $k$, the following holds trivially: $I_p \in
\mathcal{F}_q(R)$ and $\diag(\Sigma_S)=I_p$ for any $S \in\mathcal{S}$.
Write $\Sigma_S= (\sigma_{ij})$ and observe that
\[
\sum_{i\ne j} |\sigma_{ij}|^q=2ka^q.
\]
Next, we treat each case of Lemma~\ref{LEM:chi2div} separately.

\textit{Case} 1. Note first that $2ka^q=R/2<R$ so that $\Sigma_S
\in
\mathcal{F}_q(R)$. Moreover, $k^2a^4=CR^{4/q}$.

\textit{Case} 2. Note first that $2ka^q\le R/2<R$ so that $\Sigma_S
\in\mathcal{F}_q(R)$. Since $k \ge2$ and $k^2 \le R^2 n^q$, we have
\[
k \ge\frac{R}{4} \biggl(\frac{\log ({(p-1)}/{k^2}+1
)}{2n} \biggr)^{-q/2}.
\]
%
Therefore,
\[
k^2a^4\ge\frac{R^2}{16} \biggl(\frac{\log ({(p-1)}/{(R^2n^{q})}+1
)}{2n}
\biggr)^{2-q}\wedge\frac{1}4.
\]
Combining the two cases, we get
\[
k^2a^4\ge C \biggl[R^2 \biggl(
\frac{\log (
{(p-1)}/{(R^2n^{q})}+1
)}{2n} \biggr)^{2-q}\wedge R^{4/q} \wedge1 \biggr].
\]
Together with \eqref{EQ:lbtheta}, this completes the proof of the
second part of the lower bound.
\end{pf*}


\section{Extension to nonquadratic functionals}
\label{SEC:nonQuad}
Closely related to quadratic functional is the $\ell_r$ functional of
covariance matrices, which is defined by
%
\begin{equation}
\ell_r(\Sigma) = \max_{i \le p} \sum
_{j \le p} {{\llvert {{\sigma _{ij}}} \rrvert
}^r}.
\end{equation}
It is often used to measure the sparsity of a covariance matrix and
plays an important role in estimating sparse covariance matrix. This
along the theoretical interest on the difficulty of estimating such a
functional give rise to this study.
Note that $\ell_1 (\Sigma)$ functional is indeed the $\ell_1$-norm of
the covariance matrix $\Sigma$, whereas when $r=2$, $\ell_r$ functional
is the maximal row-wise quadratic functional. Thus, the nonquadratic
$\ell_r$ functional is just a natural extension of such a maximal
quadratic functional, whose optimal estimation problem will be the main
focus of this section.

\subsection{Optimal estimation of \texorpdfstring{$\ell_r$}{$ell_r$} functionals}

We consider a class of matrix with row-wise sparsity structure as follows:
%
\begin{equation}
\label{EQ:defGq} \mathcal{G}_q(R) = \biggl\{ \Sigma\in
\bS_p^+: \max_{i \le p} \sum
_{j
\le p} |\sigma_{ij}|^q \le R , \diag(
\Sigma)=I_p \biggr\} ,
\end{equation}
for $q \in[0,r)$ and $R>0$ which can depend on $n$ and $p$.
A similar class of covariance matrices has been considered by \citet
{BicLev08a} and \citet{CaiZho12}.


\begin{thmm}
\label{1595159511599}
Fix $q \in[0,r), R>0$ and assume that $2\log p < n$ and $R^2 <
(p-1)n^{-q}/2$. Then there exists a positive constant $C_4>0$ such that,
\[
\mathop{\inf} _{{{\hat L}}} \mathop{\sup} _{\Sigma\in{\mathcal
{G}_q(R)}}
\mathbb{E}_{\Sigma} \bigl[ \bigl({\hat L} - \ell_r(\Sigma )
\bigr)^2 \bigr]\ge C_4 \tilde\phi_{n,p}(q,R) ,
\]
where $\tilde\phi_{n,p}(q,R)$ is defined by
%
\begin{eqnarray}
\label{EQ:philr} \tilde\phi_{n,p}(q,R)&=& R^2
\frac{\log p}{n}
\nonumber
\\[-8pt]
\\[-8pt]
\nonumber
&&{}\boldss{\vee} \biggl\{ R^2 \biggl(
\frac{\log ({(p-1)}/{(R^2 n^{q})}+1 )}{2n} \biggr)^{r-q}\wedge R^{2r/q} \wedge1 \biggr\}
\end{eqnarray}
and the infimum is taken over all measurable functions $\hat L$ of the
sample $X_1, \ldots, X_n$.
\end{thmm}

The proof is similar to that of Theorem~\ref{TH:LBquad} and is
relegated to the \hyperref[app]{Appendix}.

As in \eqref{eq:opt}, when $1 < R^2 < p^\alpha n^{-q}$ for some
$\alpha
< 1$, the lower bound in Theorem~\ref{1595159511599} can be written as
%
\begin{equation}
\label{eq:optlr} \tilde\phi_{n,p}(q,R)= R^2
\frac{\log p}{n} \boldss{\vee} \biggl\{ R^2 \biggl(
\frac{\log p}{n} \biggr)^{r-q} \wedge1 \biggr\}.
\end{equation}

To establish the upper bound, we consider again a thresholding
estimator. Naturally, we estimate $\ell_r$ functional of each single
row, denoted by $\ell_r^{(i)}(\Sigma) = \sum_{j} {{\llvert  {{\sigma
_{ij}}} \rrvert }^r}$, using the thresholding technique. Following the
same notation as the previous section, the estimator is defined by
%
\begin{equation}
\label{Estlr} \widehat{\ell_r(\Sigma)} = \ell_r(
\tilde\Sigma_\tau)= \max_i \ell
_r^{(i)}(\tilde\Sigma_\tau) = \max
_i \sum_{j \le p} |\hat\sigma
_{ij}|^r \mathbh{1} \bigl\{ {\llvert {{{\hat
\sigma}_{ij}}} \rrvert > \tau} \bigr\} ,
\end{equation}
for a threshold $\tau> 0$. We will see in the next theorem that this
estimator achieves the adaptive minimax optimal rate.

\begin{thmm}
\label{TH:UBlr}
Assume that $\gamma\log(p)<n$ for some constant $\gamma>8$ and fix
$C_0 \ge4$. Consider the threshold
\[
\tau= 2C_0\sqrt{\frac{\gamma\log p}{n}}
\]
and assume that $\tau\le1$.
Then, for any $q \in[0,r), R>0$, the plug-in estimator $\ell_r(\tilde
\Sigma_\tau)$ satisfies
\[
\mathbb{E} \bigl[ \bigl(\ell_r(\tilde\Sigma_\tau) -
\ell_r(\Sigma ) \bigr)^2 \bigr] \le C_5
\tilde\psi_{n,p}(q,R) + C_6 p^{4-\gamma/2} ,
\]
where
\[
\tilde\psi_{n,p}(q,R)= \cases{ %
\displaystyle\frac{R^{2} \log p}{n},&\quad $\mbox{if } q < \max\{r-1,0\},$
\vspace*{2pt}\cr
\displaystyle R^2 \biggl(\frac{\gamma\log p}{n} \biggr)^{r-q}, & \quad
$\mbox{if } q
\ge\max \{ r-1,0\}$}
\]
and $C_5$ and $C_6$ are positive constants.
\end{thmm}

The proof of this theorem is a generalization of the proof of Theorem~\ref{TH:UBquad} but some aspects that have independent value are
presented here. In the proof of Theorem~\ref{TH:UBquad}, we used the
decomposition
\[
\hat\sigma_{ij}^2- \sigma_{ij}^2 = 2
\sigma_{ij}(\hat\sigma _{ij}-\sigma _{ij})+ (\hat
\sigma_{ij}-\sigma_{ij})^2 ,
\]
which is actually the Taylor expansion of $\hat\sigma_{i,j}^2$ at
$\sigma_{i,j}$. Carefully scrutinizing the proof, we find that the
first term has the parametric rate $O(R^2/n)$ whereas the second term
contributes to the rate $O(R^2(\log p/n)^{2-q})$. This phenomenon can
be generalized to the $\ell_r$-functional. In the latter case, we will
apply the Taylor expansion of $|\hat\sigma_{ij}|^r$ at $|\sigma_{ij}|$, and
the first-order term will contribute to the parametric rate of
$O(R^2\log p/n)$ while the second-order term has the rate $O(R^2(\log
p/n)^{r-q})$. The elbow effect stems from the dominance of estimation
errors of the first- and second-order terms of Taylor's expansion. We
relegate the complete proof to the supplementary material.

A few remarks should be mentioned:
\begin{longlist}[1.]
\item[1.] The combination of the two theorems imply that the estimator
$\ell
_r(\tilde\Sigma_\tau)$ is minimax adaptive over the space $\{
\mathcal
{G}_q(R), q\in[0,r), R>0\}$ under very mild conditions. The adaptive
minimax optimal rate of convergence is given by \eqref{eq:optlr}. The
term $p^{4-\gamma/2}$ can be made arbitrarily small by choosing large
enough $\gamma$.

\item[2.] The $\ell_r$ functional involves the maxima of the row sums.
Compared it with estimating the quadratic functional, we need to pay
the price of an extra $\log p$ term in the parametric rate.

\item[3.] The rate $\tilde\phi_{n,p}(q,R)$ presents the elbow phenomenon
at $q=r-1$ if $r > 1$. So quadratic row-wise functional $\ell_2(\tilde
\Sigma_\tau)$ bears the same elbow behavior as the quadratic functional
$Q(\tilde\Sigma_\tau)$. 
\end{longlist}

\subsection{Optimal detection of correlations} In this subsection, we
illustrate the intrinsic link between functional estimation and
hypothesis testing. To that end, consider the following hypothesis
testing problem:
\begin{eqnarray*}
&& H_0:\qquad X \sim\mathcal N(0, I_p) ,
\\
&& H_1:\qquad X \sim\mathcal N\bigl(0, I_p + \kappa\cdot
\offdiag(\Sigma)\bigr) , \qquad\Sigma\in\bigcup_{q \in[0,r)} \bigl
\{\mathcal{G}_q(R): \ell _r\bigl(\offdiag (\Sigma)\bigr)
= 1 \bigr\}.
\end{eqnarray*}
This problem is intimately linked to sparse principal component
analysis [\citeauthor{BerRig13a} (\citeyear{BerRig13a,BerRig13b})]. A natural question associated
with this problem is to find the minimal signal strength $\kappa$ such
that these hypotheses can be tested with high accuracy.

The previous subsection provides the optimal estimate for $\ell
_r(\offdiag(\Sigma))$. However, we need a result with high probability
rather than in expectation. Using Lemma~4.2 in the supplementary material [\citet
{FanRigWan15}] and arguments similar to those employed to prove
Theorem~\ref{TH:UBlr}, it is not hard to show that
\[
\bigl| \ell_r(\tilde\Sigma_\tau) - \ell_r(
\Sigma) \bigr| \le CR \biggl(\frac
{\gamma\log p}{n} \biggr)^{{(r-q)}/{2}} = CR \biggl(
\frac{2\log p +
\log
(4/\delta)}{n} \biggr)^{{(r-q)}/{2}} ,
\]
with probability larger than $1-4p^{-(\gamma-2)} =: 1-\delta$.
Therefore, letting
\begin{eqnarray*}
s_0 &=& 1 + CR \biggl(\frac{2\log p + \log(4/\delta)}{n} \biggr)^{{(r-q)}/{2}} ,
\\
s_1 &=& 1+\kappa^r - CR \biggl(\frac{2\log p + \log(4/\delta
)}{n}
\biggr)^{{(r-q)}/{2}} ,
\end{eqnarray*}
we get $\mathbb P_{H_0} ( \ell_r(\tilde\Sigma_\tau) \le s_0) \ge1 -
\delta$ and $\mathbb P_{H_1} (\ell_r(\tilde\Sigma_\tau) \ge s_1)
\ge1
- \delta$. Here, $\mathbb P_{H_0}$ denotes the probability under the
null hypothesis and $\mathbb P_{H_1}$ denotes the largest probability
over the composite alternative. To build a hypothesis test, note that
if $s_1 > s_0$, then for any $s \in[s_0, s_1]$, the test $\psi
=\mathbh
{1}\{\ell_r(\tilde\Sigma_\tau) \ge s\}$ satisfies $\mathbb
P_{H_0}(\psi
=1)\vee\mathbb P_{H_1}(\psi=0) \le\delta$. We say that the test
$\psi
$ discriminates between $H_0$ and $H_1$ with accuracy $\delta$.

\begin{thmm}
\label{TH:UBdetection}
Assume that $n, p, R, q, r$ and $\delta$ are such that $\bar\kappa<
1$ where
\[
\bar\kappa:= 2CR^{1/r} \biggl(\frac{2\log p + \log(4/\delta)}{n} \biggr)^{{(r-q)}/{(2r)}}.
\]
Then, for any $\kappa> \bar\kappa$ and for any $s \in[s_0, s_1]$, the
test $\psi=\mathbh{1}\{\ell_r(\tilde\Sigma_\tau) \ge s\}$
discriminates between $H_0$ and $H_1$ with accuracy $\delta$.
\end{thmm}

\begin{figure}[b]

\includegraphics{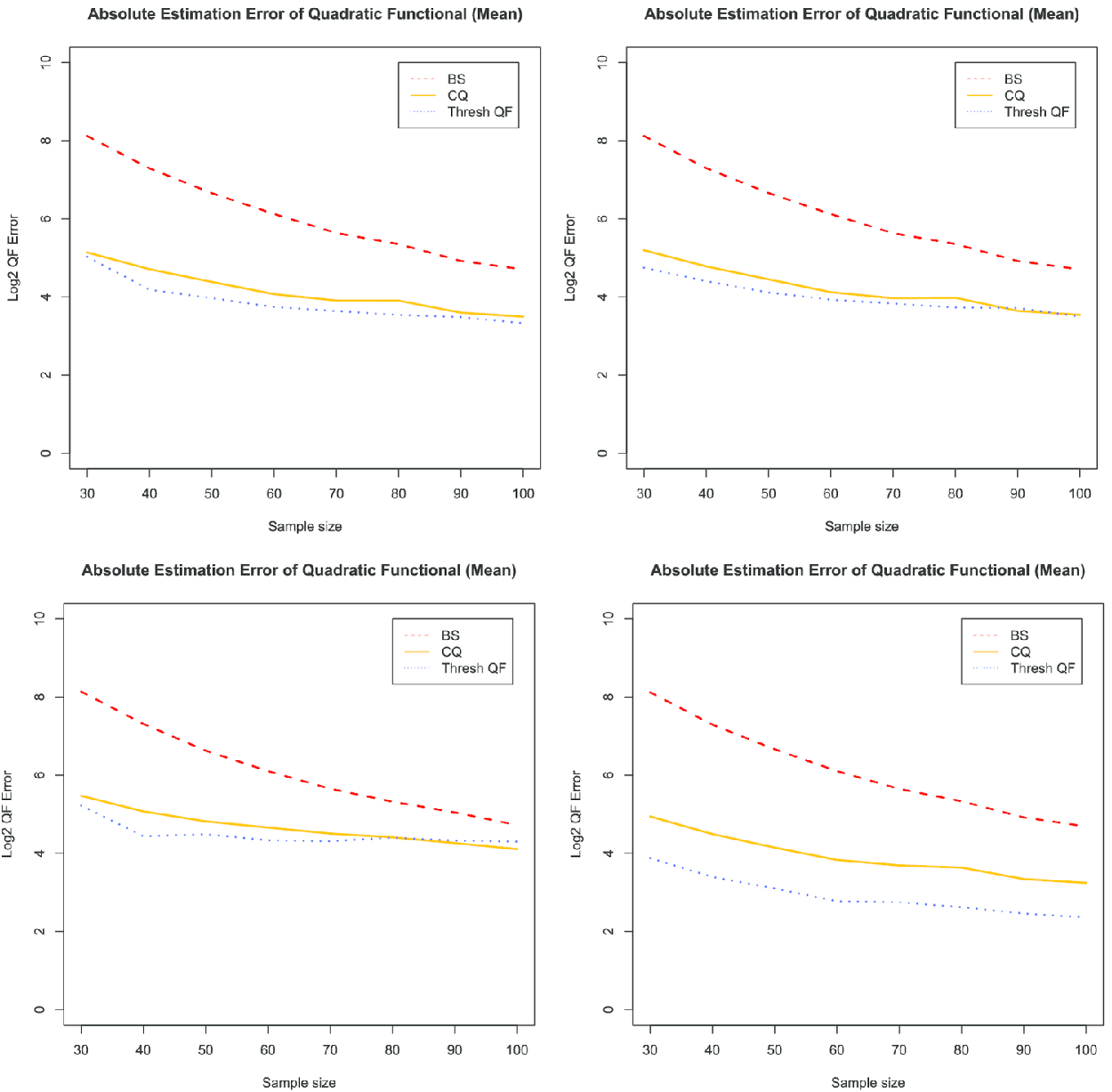}

\caption{Performance of estimating $\|\Sigma\|_\mathsf{F}^2$ using
thresholded
estimator $\hat Q + \hat D$ (dotted), CQ (solid) and BS
(dashed). The mean of absolute errors over 500 repetitions {in log
scale (base 2)} versus the sample size were reported for matrix \textup{M1} (top
left), \textup{M2} (top right), \textup{M3} (bottom left), \textup{M4} (bottom right).}
\label{fig:behavior}
\end{figure}

The minimax risk for the correlation detection is given in the next
theorem, which will be proved in the \hyperref[app]{Appendix}.

\begin{thmm}
\label{35753575357}
For fixed $\nu> 0$, define $\underline\kappa> 0$ by
\[
\underline\kappa:= R^{1/r} \biggl( \frac{\log(\nu
p/(R^2n^q))}{2n}
\biggr)^{{(r-q)}/{(2r)}}.
\]
Then for any $\kappa< \underline\kappa$,
\[
\inf_{\psi} \bigl\{ \mathbb P_{H_0} (\psi= 1) \vee
\mathbb P_{H_1} (\psi= 0) \bigr\} \ge C_{\nu} ,
\]
where the infimum is taken over all possible tests and $C_{\nu}>0$ is a
continuous function of $\nu$ that tends to $1/2$ as $\nu\to0$.
\end{thmm}

If we assume the high-dimensional regime $R^2 < p^\alpha n^{-q}$ for
some $\alpha< 1$ as discussed before, then the lower bound matches the
upper bound. So the theorem concludes that no test has asymptotic power
for correlation detection unless $\kappa$ is of higher order than
$R^{1/r} (\log p/n)^{{(r-q)}/{(2r)}}$ and the detection method based on
optimal $\ell_r(\Sigma)$ estimation is also optimal for testing
existence of correlation.

\section{Numerical experiments}
\label{SEC:num}
Simulations are conducted in this section to evaluate the numerical
performance of our plug-in estimator for quadratic functionals. Then
the proposed method is applied to two high-dimensional testing
problems: simulated two-sample data and real financial equity market data.

\subsection{Quadratic functional estimation}
\begin{figure}

\includegraphics{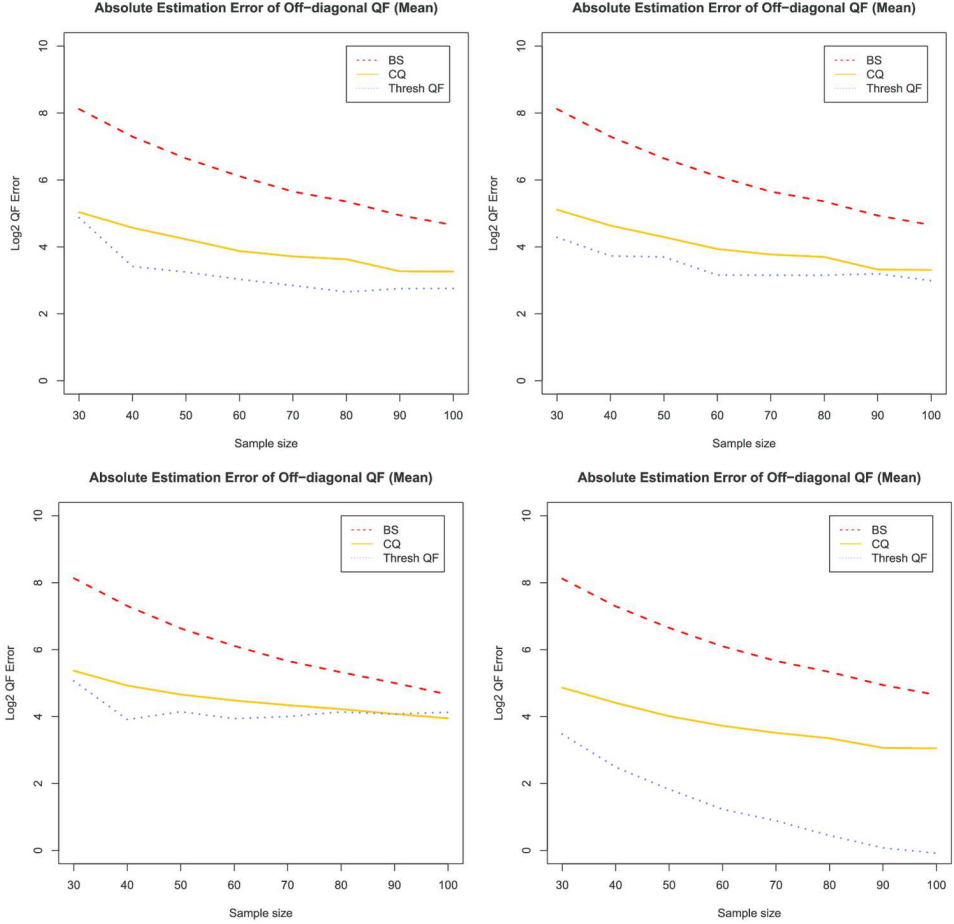}

\caption{Performance of estimating $Q(\Sigma)$ using thresholded
estimator $\hat Q$ (dotted), CQ-$\hat D$ (solid) and
BS-$\hat D$ (dashed). The mean of absolute errors over 500
repetitions {in log scale (base 2)} versus the sample size were
reported for matrix \textup{M1} (top left), \textup{M2} (top right),
\textup{M3} (bottom left), \textup{M4}
(bottom right).}
\label{fig:behavior_Q}
\end{figure}

We first study the behavior of estimators $\widehat{Q(\Sigma)} +
\widehat{D(\Sigma)}$ for the total quadratic functional and $\widehat
{Q(\Sigma)} = Q(\tilde\Sigma_\tau)$ for its off-diagonal part. To that
end, four sparse covariance matrix structures were used in the simulations:
\begin{itemize}
\item[(M1)] auto-correlation AR($1$) covariance matrix $\sigma_{ij} =
0.25^{\llvert i-j\rrvert }$;
\item[(M2)] banded correlation matrix with $\sigma_{ij} = 0.3$ if
$\llvert i-j\rrvert =1$ and $0$ otherwise;
\item[(M3)] sparse matrix with a block, size $p/20$ by $p/20$, of
correlation $0.3$;
\item[(M4)] identity matrix (it attains the maximal level of sparsity).
\end{itemize}
We chose $p = 500$ and let $n$ vary from $30$ to $100$. For estimating
the total quadratic functional, our proposed thresholding estimator, BS
[\citet{{BaiSar96}}] estimator and CQ [\citet{CheQin10}] estimator were
applied to each setting for repetition of $500$ times. Their mean
absolute estimation errors were reported {in log scale (base 2)} in
Figure~\ref{fig:behavior} with their standard deviations omitted here.
BS and CQ cannot be directly used for off-diagonal quadratic functional
estimation, so we deducted $\widehat{D(\Sigma)}$ from both of them to
serve as an estimator for only the off-diagonal part. The mean absolute
estimation errors, compared with our proposed estimator $Q(\tilde
\Sigma
_\tau)$, are depicted {in log scale (base 2)} in Figure~\ref{fig:behavior_Q}.

%

%

The four plots correspond to the aforementioned four covariance
structures. We did not report the estimation error of directly using
the naive plug-in which is an obvious disaster. In all the four cases,
the BS (dashed line) method does not perform well in the ``large
$p$ small $n$'' regime. The method CQ (solid line) exhibits a
relatively small estimation error in general, but it can still be
improved using the thresholding method. As theory shows, the method CQ
is ratio-consistent [\citet{CheQin10}], so our method (dotted line)
is better only to a second order, which was captured by the small gap
between dotted and solid curves. {When estimating only off-diagonal
quadratic functionals (Figure~\ref{fig:behavior_Q}), the advantage of
the thresholding method is even sharper since the error caused by
nonsparse diagonal elements has been eliminated.}
The improved performance comes from the prior knowledge of sparsity,
thus our method works best for very sparse matrix, especially well for
identity matrix as seen in Figure~\ref{fig:behavior}.

A practical question is how to choose a proper threshold, as this is
important to the performance of the thresholding estimator. In the
above simulations, we chose $\tau= C\sqrt{\log p/n}$ with constant $C$
slightly different for the four cases but all close to $1.5$. In the
next two applications to hypothesis testing, we employ the cross
validation to choose a proper thresholding. The procedure consists of
the following steps:
\begin{longlist}[(1)]
\item[(1)] The data is split into training data $D_{\mathcal S}^{(v)}$
of sample size $n_1$ and testing data $D_{\mathcal S^c}^{(v)}$ of
sample size $n-n_1$ for $m$ times, $v=1,2,\ldots,m$.
\item[(2)] The training data $D_{\mathcal S}^{(v)}$ is used to
construct the thresholding estimator $Q(\tilde\Sigma_\tau^{(v)})$
under a sequence of thresholds while the testing data $D_{\mathcal
S^c}^{(v)}$ constructs the nonthresholded ratio-consistent estimator
$\hat Q^{(v)}$, for example, using CQ estimator of ${\|\Sigma\|
_\mathsf{F}^2}$.
\item[(3)] The candidates of thresholds are $\tau_j = j \Delta\sqrt
{\log(p)/n_1}$ for $j=1,2,\ldots,J$ where $J$ is chosen to be a
reasonably large number, say $50$, and $\Delta$ is such that $J \Delta
\sqrt{\log(p)/n_1} \le\hat M:= \max_{i} \hat\sigma_{ii}$.
\item[(4)] The final $j^*$ is taken to be the minimizer of the
following problem:
\[
\min_{j \in \{1,2,\ldots,J \} } \frac{1}{m} \sum
_{v = 1}^{m} \bigl\llvert Q\bigl(\tilde
\Sigma_{\tau_j}^{(v)}\bigr) - \hat Q^{(v)} \bigr\rrvert .
\]
\item[(5)] The final estimator $Q(\tilde\Sigma_{\tau_{j^*}})$ is
obtained by applying threshold $\tau_{j^*} = j^* \Delta\sqrt{\log
(p)/n}$ to the empirical covariance matrix of the entire $n$ data.
\end{longlist}
\citet{BicLev08a} suggested to use $n_1 = n/\log{n}$ for estimating
covariance matrices. This is consistent with our experience for
estimating functionals when no prior knowledge about the covariance
matrix structure is provided. We will apply this splitting rule in the
later simulation studies on high-dimensional hypothesis testing.

\subsection{Application to high-dimensional two-sample testing}
\label{SUBSEC:testing}
In this section, we apply the thresholding estimator of quadratic
functionals to the high-dimensional two-sample testing problem. Two
groups of data are simulated from the Gaussian models:
\[
X_{i,j} \sim\mathcal{N}(\mu_i, \Sigma)\qquad \mbox{for $i=1,2$
and $j=1,\ldots, n/2$}.
\]
The dimensions considered for this problem are $(p,n)\in \{(500,100),\break
(1000,150),   (2000,200)\}$. For simplicity, we choose $\Sigma$ to be a
correlation matrix and choose the sparse covariance structure to be $2$
by $2$ block diagonal matrices with $250$ of them having correlations
$0.3$ and the rest having correlations $0$. So the off-diagonal
quadratic functional is always $45$, which does not increase with $p$
in our setting. The mean vectors $\mu_1$ and $\mu_2$ are chosen as
follows. Let $\mu_1 = 0$ and the percentage of $\mu_{1,k} = \mu_{2,k}$
to be in $\{0\%, 50\%, 95\%, 100\%\}$. The $100\%$ proportion
corresponds to the case where the two groups are identical, thus gives
information about accuracy of the size of the tests. The $95\%$
proportion represents the situation where the alternative hypotheses
are sparse. For those $k$ such that $\mu_{1,k} \ne\mu_{2,k}$, we
simply chose the value of each $\mu_{2,k}$ equally. To make the power
comparable among different configurations, we use a constant
signal-to-noise ratio $\eta= \llVert \mu_1-\mu_2\rrVert /\sqrt{\tr
(\Sigma^2)} = 0.1$ across experiments.

\begin{table}[b]
\caption{Empirical testing power and size of $6$ testing methods based
on $500$ simulations}
\label{tab:TestingSim}
\begin{tabular*}{\textwidth}{@{\extracolsep{\fill}}lcccccc@{}}
\hline
\textbf{Prop. of equalities} & \textbf{BS} & \textbf{newBS} &
\textbf{CQ} & \textbf{newCQ} & \textbf{Bonf} & \textbf{BH} \\
\hline
& \multicolumn{6}{c}{$p=500, n=100$} \\
0\% & 0.408 & 0.422 & 0.428 & 0.432 & 0.104 & 0.110 \\
50\% & 0.396 & 0.422 & 0.418 & 0.428 & 0.110 & 0.116 \\
95\% & 0.422 & 0.440 & 0.438 & 0.442 & 0.208 & 0.214 \\
100\% (size) & 0.030 & 0.036 & 0.036 & 0.038 & 0.042 & 0.042 \\[3pt]
& \multicolumn{6}{c}{$p=1000, n=150$} \\
0\% & 0.696 & 0.710 & 0.718 & 0.718 & 0.082 & 0.086 \\
50\% & 0.698 & 0.712 & 0.712 & 0.714 & 0.106 & 0.112 \\
95\% & 0.702 & 0.716 & 0.718 & 0.722 & 0.308 & 0.328 \\
100\% (size) & 0.040 & 0.044 & 0.048 & 0.046 & 0.050 & 0.050 \\[3pt]
& \multicolumn{6}{c}{$p=2000, n=200$} \\
0\% & 0.930 & 0.938 & 0.940 & 0.940 & 0.138 & 0.146 \\
50\% & 0.918 & 0.922 & 0.924 & 0.928 & 0.104 & 0.106 \\
95\% & 0.922 & 0.928 & 0.930 & 0.930 & 0.324 & 0.338\\
100\% (size) & 0.046 & 0.050 & 0.050 & 0.050 & 0.046 & 0.046 \\
\hline
\end{tabular*}
\end{table}

Table~\ref{tab:TestingSim} reports the empirical power and size of six
testing methods based on $500$ repetitions.


%
\begin{longlist}[\quad]
\item[(BS)]  Bai and Saranadasa's original test.

\item[(newBS)] Bai and Saranadasa's modified test where $\tr(\Sigma^2)$ is
estimated by thresholding the sample covariance matrix.

\item[(CQ)] Chen and Qin's original test.

\item[(newCQ)] Chen and Qin's modified test where $\tr(\Sigma_i^2)$ and
$\tr
(\Sigma_1 \Sigma_2)$ are estimated by thresholding their empirical
counterparts.

\item[(Bonf)] Bonferroni correction: This method regards the high-dimensional
testing problem as $p$ univariate testing problems. If there is a
$p$-value that is less than $0.05/p$, the null hypothesis is rejected.

\item[(BH)] Benjamini--Hochberg method. The method is similar to the
Bonferroni correction, but employs the Benjamini--Hochberg method in
decision making.
\end{longlist}

For estimating quadratic functionals, the cross-validation is employed
using $n/\log(n)$ splitting rule. The first four methods are evaluated
at the $5\%$ significance level while Bonferroni correction and
Benjamini--Hochberg correction are evaluated at $5\%$ family-wise error
rate or FDR. We also list the average relative estimation errors for
the quadratic functionals of the first four methods in Table~\ref
{tab:TestingErr}. Here, the average is taken over four different
proportions of equalities and the average for CQ and newCQ is also
taken over errors in estimating $\tr(\Sigma_1^2)$ and $\tr(\Sigma_2^2)$.

\begin{table}
\caption{Mean and SD of relative errors for estimating quadratic
functionals (in percentage)}
\label{tab:TestingErr}
\begin{tabular*}{\textwidth}{@{\extracolsep{\fill}}lccc@{}}
\hline
& \multicolumn{1}{c}{$\bolds{p=500, n=100}$} &
\multicolumn{1}{c}{$\bolds{p=1000, n=150}$} &
\multicolumn{1}{c@{}}{$\bolds{p=2000, n=200}$} \\
\hline
BS & 4.93 (2.48) & 4.47 (1.56) & 5.05 (1.10) \\
newBS & 2.12 (1.43) & 0.74 (0.56) & 0.54 (0.40) \\
CQ & 3.72 (1.97) & 2.32 (1.24) & 1.70 (0.91) \\
newCQ & 2.77 (1.38) & 1.27 (0.64) & 0.62 (0.33) \\
\hline
\end{tabular*}
\end{table}

Several comments are in order. First, the first four methods based on
Wald-type of statistic with correlation ignored perform much better, in
terms of the power, than the last two methods which combines individual
tests. Even in the case that proportional of equalities is $95\%$ where
the individual difference is large for nonidentical means, aggregating
the signals together in the Wald-type of statistic still outperforms.
However, in the case of $0\%$ identical means, the power of Bonferroni
or FDR method is extremely small, due to small individual differences.
Second, the method newCQ, which combines CQ and thresholding estimator
of the quadratic functional, has the highest power and performs the
best among all methods. The corrected BS method also improves the
performance by estimating the quadratic functionals better compared
with original BS. CQ indeed is more powerful than BS as claimed by
\citet
{CheQin10}, but we can even improve the performance of those two
methods more by leveraging the sparsity structure of covariance matrices.

\subsection{Estimation of \texorpdfstring{$\ell_r$}{$ell_r$} functional and correlation detection}

In order to check the effectiveness of using $\ell_r$ norm of the
thresholded sample matrix to detect correlation, let us take one simple
matrix structure as an example and use $r=1$. Under $H_0$, assume $X
\sim\mathcal N(0, I_p)$; while under $H_1$, $X \sim\mathcal
N(0,\Sigma
)$, where $\Sigma_{ij} = 0.8$ if $i,j \in\mathcal S$ and $\mathcal S$
is a random subset of size $p/20$ in $\{1, 2,\ldots , n\}$. We chose to
use $p = 500$ and generated $n = 100$ independent random vectors under
both $H_0$ and $H_1$. The whole simulation was done for $N = 1000$ times.

\begin{figure}

\includegraphics{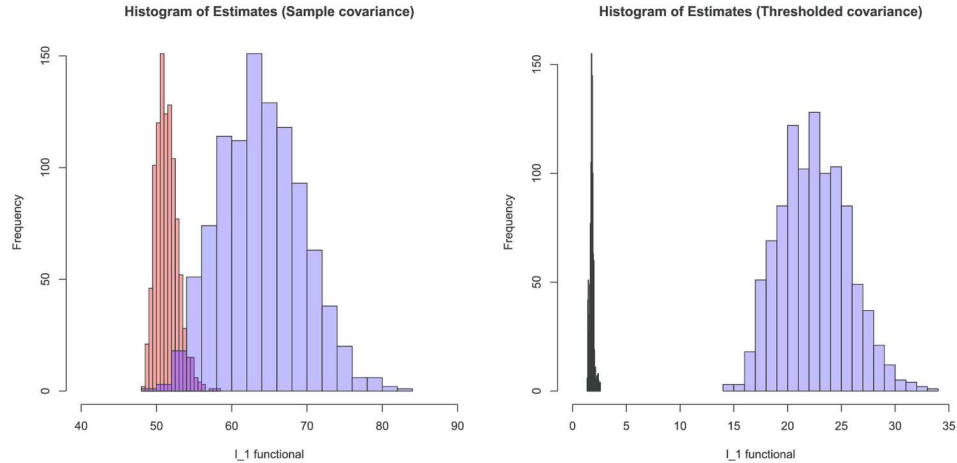}

\caption{Histogram of $1000$ $\ell_1$ functional estimates for $H_0$
 and $H_1$  by $\ell_r(\hat\Sigma)$ (left) and optimal
estimator $\ell_r(\tilde\Sigma_{\tau})$ (right).}
\label{fig:detection}
\end{figure}

We compare the $\ell_1$ norm estimates based on empirical covariance
matrix $\ell_1(\hat\Sigma)$ and thresholded empirical covariance matrix
$\ell_1(\tilde\Sigma_{\tau})$. The threshold is decided by cross
validation with $n/\log(n)$ splitting. The simulations yielded $N$
estimates for both null and alternative hypotheses, which were plotted
in Figure~\ref{fig:detection}. The optimal estimator $\ell_1(\tilde
\Sigma_{\tau})$ perfectly discriminates the null and alternative
hypotheses while $\ell_1(\hat\Sigma)$ overestimates $\ell_1$ functional
and blurs the difference of the two hypotheses.

\subsection{Application to testing multifactor pricing model}

\begin{figure}

\includegraphics{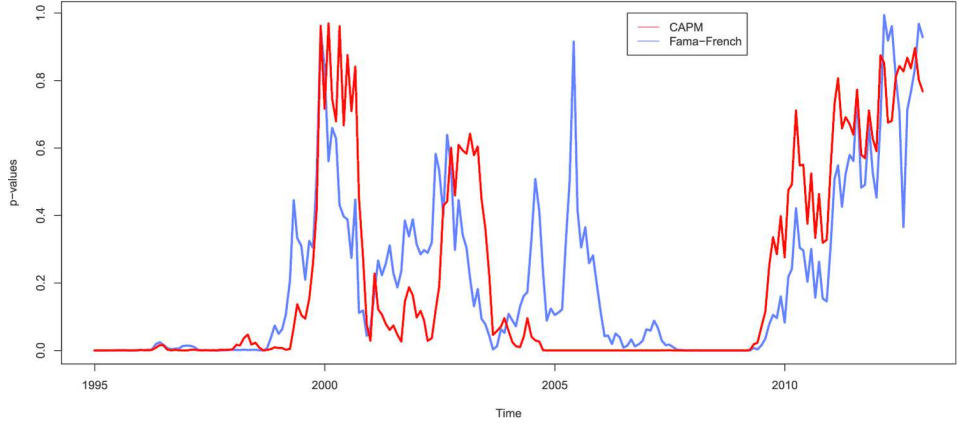}

\caption{$P$-values of testing $H_0: \alpha= 0$ in the CAMP and
Fama--French 3 factor models based on S\&P 500 monthly returns from
January 1995 to December 2012.}
\label{fig:sp500}
\end{figure}

In this section, we test the validity of the Capital Asset Pricing
Model (CAPM) and Fama--French models using Pesaran and Yamagata's method
\eqref{eq2.3} for the securities in the Standard \& Poor 500 (S\&P 500)
index. Following the literature, we used 60 monthly stock returns to
construct test statistics since monthly returns are nearly independent.
The composition of index keeps changing annually, so we selected only
$276$ large stocks. The monthly returns (adjusted for dividend) between
January 1990 and December 2012 are downloaded from the Wharton Research
Data Services (WRDS) database. The time series on the risk-free rates
and Fama--French three factors are obtained from Ken French's data
library. If only the first factor, that is, the excessive return of the
market portfolio is used, the Fama--French model reduces to the CAPM
model. We tested the null hypothesis $H_0: \alpha= 0$ for both models.
The $p$-values of the tests are depicted in Figure~\ref{fig:sp500}, which
are computed based on running windows of previous 60 months.

The results suggest that market efficiency is time dependent and the
Fama--French model are rejected less frequently than the CAPM. Before
1998, the evidence that $\alpha\ne0$ is very strong. After 1998, the
Fama--French 3-factor model holds most of the time except the period
2007--2009 that contains the financial crisis. On the other hand, the
CAPM is rejected for an extended period of time during this period.

\begin{appendix}\label{app}
\section{A technical lemma on \texorpdfstring{$\chi^2$}{$chi^2$} divergences}
\label{SUB:proof:LEM:chi2div}

\begin{lem}
\label{LEM:chi2div}
Consider a mixture of Gaussian product distributions
\[
\tilde{\mathbb{P}}^n=\frac{1}{m}\sum
_{j=1}^m \mathbb{P}_j^n,
\]
where $\mathbb{P}_j \sim\mathcal{N}_p(0, \Sigma_j)$ such that
$\mathbb
{P}_j \ll\mathbb{P}_0$. Then
%
\begin{equation}
\label{EQ:LEM:chi2-1} \chi^2\bigl(\tilde{\mathbb{P}}^n,
\mathbb{P}_0^n\bigr)=\frac{1}{m^2} \sum
_{j,k=1}^m\bigl|I_p-(\Sigma_j-I)
(\Sigma_k-I)\bigr|^{-n/2}-1.
\end{equation}
Furthermore, assume $2(\log p) \le n$. Consider the mixture $\bar
{\mathbb{P}}^n$ defined in the proof of Theorem~\ref{TH:LBquad} where
$k$ and $a$ are defined as follows:

1. If $R<4(\frac{\log p}{n})^{q/2}$, then take $k=1$ and
$a=(R/4)^{1/q}$.

2. If $R\ge4(\frac{\log p}{n})^{q/2}$, then take $k$ to be
the largest integer such that
%
\begin{equation}
\label{EQ:choice_k} k\le\frac{R}{2} \biggl(\frac{\log ({(p-1)}/{k^2}+1 )}{2n}
\biggr)^{-q/2}
\end{equation}
and
%
\begin{equation}
\label{EQ:choice_a} a= \biggl(\frac{\log ({(p-1)}/{k^2}+1 )}{2n} \biggr)^{1/2}\wedge
(2k)^{-1/2}.
\end{equation}
Such choices yield in both cases
%
\begin{equation}
\label{EQ:LEM:chi2-2} \chi^2\bigl(\bar{\mathbb{P}}^n,
\mathbb{P}_0^n\bigr)\le e-1.
\end{equation}
Moreover, in case 2 we have that \textup{(i)} $k \ge2$ and \textup{(ii)} under the
assumption that $R^2 < (p-1) n^{-q} / 2$, we also have $k^2 \le R^2 n^q
< (p-1)/2$.
\end{lem}

\begin{pf}
To unify the notation, we will work directly with $\mathbb{P}_S, S\in
\mathcal{S}$ rather than $\mathbb{P}_j, j \in[m]$. However, in the
first part of the proof, we will not use the specific form $\Sigma_S$
nor that of $\mathcal{S}$. For now, we simply assume that $\Sigma_S$ is
invertible (we will check this later on). Recall that
\[
\chi^2\bigl(\bar{\mathbb{P}}^n, \mathbb{P}_0^n
\bigr)=\mathbb{E}_0 \biggl[ \biggl(\frac
{\mathrm{d} \bar{\mathbb{P}}^n}{\mathrm{d} \mathbb{P}_0^n}-1
\biggr)^2 \biggr]=\frac{1}{|\mathcal{S}|^2}\sum_{S,T \in\mathcal{S}}
\biggl(\mathbb {E}_0 \biggl[\frac{\mathrm{d} \mathbb{P}_S}{\mathrm{d} \mathbb
{P}_0}\frac
{\mathrm{d} \mathbb{P}_T}{\mathrm{d} \mathbb{P}_0}
\biggr] \biggr)^n-1 ,
\]
where $\mathbb{E}_0$ denotes the expectation with respect to $\mathbb
{P}_0$. Furthermore,
\[
\mathbb{E}_0 \biggl[\frac{\mathrm{d} \mathbb{P}_S}{\mathrm{d}
\mathbb
{P}_0}\frac{\mathrm{d} \mathbb{P}_T}{\mathrm{d} \mathbb{P}_0} \biggr]=
\frac
{1}{(|\Sigma_S||\Sigma_T|)^{1/2}}\mathbb{E}_0 \biggl[\exp \biggl(-
\frac{1}2X^\top\bigl(\Sigma_S^{-1}+
\Sigma_T^{-1}-2I_p \bigr)X \biggr) \biggr].
\]
Consider the spectral decomposition of $\Sigma_S^{-1}+\Sigma
_T^{-1}-2I_p=U\Lambda U^\top$, where $U$ is an orthogonal matrix and
$\Lambda$ is a diagonal matrix with eigenvalues $\lambda_1, \ldots,
\lambda_p$ on its diagonal. Then, by rotational invariance of the
Gaussian distribution, it holds
\begin{eqnarray*}
&&\mathbb{E}_0 \biggl[\exp \biggl(-\frac{1}2
X^\top\bigl(\Sigma_S^{-1}+\Sigma
_T^{-1}-2I_p \bigr)X \biggr) \biggr]\\
&&\qquad=
\mathbb{E}_0 \biggl[\exp \biggl(-\frac{1}2 X^\top
\Lambda X \biggr) \biggr]
\\
&&\qquad=\prod_{j=1}^p\mathbb{E}_0
\biggl[\exp \biggl(-\frac{1}2 \lambda_j X_j^2
\biggr) \biggr]
\\
&&\qquad=\cases{ %
\displaystyle\prod_{j=1}^p(1+
\lambda_j)^{-1/2}=|I+\Lambda|^{-1/2},& \quad$\mbox{if } \max
_j \lambda_j <1,$
\vspace*{2pt}\cr
\infty,& \quad$\mbox{otherwise}.$}
\end{eqnarray*}
To ensure that the above expression is finite, note that the
Cauchy--Schwarz inequality yields
\begin{eqnarray*}
\biggl(\mathbb{E}_0 \biggl[\frac{\mathrm{d} \mathbb{P}_S}{\mathrm{d}
\mathbb
{P}_0}\frac{\mathrm{d} \mathbb{P}_T}{\mathrm{d} \mathbb{P}_0}
\biggr] \biggr)^2 & \le& \mathbb{E}_0 \biggl[ \biggl(
\frac{\mathrm{d} \mathbb
{P}_S}{\mathrm
{d} \mathbb{P}_0} \biggr)^2 \biggr]\mathbb{E}_0 \biggl[
\biggl(\frac{\mathrm{d}
\mathbb{P}_T}{\mathrm{d} \mathbb{P}_0} \biggr)^2 \biggr]
\\
& = & \bigl(\chi^2(\mathbb{P}_S, \mathbb{P}_0)+1
\bigr) \bigl(\chi^2(\mathbb{P}_T, \mathbb{P}_0)+1
\bigr)<\infty,
\end{eqnarray*}
where the two $\chi^2$ divergences are finite because $\mathbb{P}_S
\ll
\mathbb{P}_0$ for any $S \in\mathcal{S}$. Therefore,
\[
\mathbb{E}_0 \biggl[\frac{\mathrm{d} \mathbb{P}_S}{\mathrm{d}
\mathbb
{P}_0}\frac{\mathrm{d} \mathbb{P}_T}{\mathrm{d} \mathbb{P}_0} \biggr]=
\frac
{|I+\Lambda|^{-1/2}}{(|\Sigma_S||\Sigma_T|)^{1/2}}=\frac
{|\Sigma
_S^{-1}+\Sigma_T^{-1}-I_p|^{-1/2}}{(|\Sigma_S||\Sigma_T|)^{1/2}}.
\]
Next, observe that
\begin{eqnarray*}
\bigl(|\Sigma_S||\Sigma_T|\bigl|\Sigma_S^{-1}+
\Sigma _T^{-1}-I_p\bigr| \bigr)^{-1/2}&=&
\bigl(\bigl(\bigl|\Sigma_S\bigl(\Sigma_S^{-1}+
\Sigma_T^{-1}-I_p\bigr)\Sigma _T
\bigr)\bigr| \bigr)^{-1/2}
\\
&=&\bigl|I-(\Sigma_S-I) (\Sigma_T-I)\bigr|^{-1/2}.
\end{eqnarray*}
%
Since we have not used the specific form of $\Sigma_S$, $S \in
\mathcal
{S}$, this bound is valid for any mixture and completes the proof
of \eqref{EQ:LEM:chi2-1}.

Next, we apply this bound to the specific choice for $\Sigma_S$
of \eqref{EQ:defSigmaS}. Note that the minimal eigenvalue of the
matrices $\Sigma_S, S \in\mathcal{S}$ is $1-\sqrt{ka^2}$. Later we
will show $2a^2 k\le1$, which implies that $\Sigma_S$ is always
positive definite. In particular, this implies that $\mathbb{P}_S \ll
\mathbb{P}_0$ for any $S \in\mathcal{S}$. Moreover, it follows from
definition \eqref{EQ:defSigmaS} that
\[
I-(\Sigma_S-I) (\Sigma_T-I)=\pmatrix{1-a^2 \bone_S^\top
\bone_T &0
\vspace*{2pt}\cr
0 & I-a^2 \bone _S\bone _T^\top} ,
\]
where $0$ is a generic symbol to indicate space filled by zeros.
Expanding the determinant along the first row (or column), we get
\begin{eqnarray*}
\bigl|I-(\Sigma_S-I) (\Sigma_T-I)\bigr|^{-1/2}&=&
\bigl(1-a^2 \bone_S^\top\bone _T
\bigr)^{-1/2}\bigl| I-a^2 \bone_S\bone_T^\top\bigr|^{-1/2}\\
&=&
\bigl(1-a^2 \bone_S^\top \bone _T
\bigr)^{-1} ,
\end{eqnarray*}
where in the second equality, we used Sylvester's determinant theorem.
By \eqref{EQ:LEM:chi2-1}, we have
\[
\chi^2\bigl(\bar{\mathbb{P}}^n, \mathbb{P}_0^n
\bigr)=\frac{1}{|\mathcal
{S}|^2}\sum_{S,T \in\mathcal{S}}
\bigl(1-a^2 \bone_S^\top\bone_T
\bigr)^{-n}-1. 
\]
As to be verified later, $2a^2k\le1$.
Using the fact that $(1-x)^{-1} \leq\exp(2x)$ for $x \in[0,1/2]$ and
the symmetry, we have
\[
\chi^2\bigl(\bar{\mathbb{P}}^n, \mathbb{P}_0^n
\bigr) \leq\frac{1}{|\mathcal
{S}|}\sum_{S \in\mathcal{S}}\exp
\bigl(2na^2 \bigl|S \cap [k] \bigr| \bigr) - 1 =\mathcal{E}\bigl[\exp
\bigl(2na^2 \bigl|S \cap [k]\bigr | \bigr) -1 \bigr],
\]
where $\mathcal{E}$ denotes the expectation with respect to the
distribution of $S$ randomly chosen from $\mathcal{S}$. In particular,
$|S\cap[k]|=\sum_{i=1}^{k}\mathbh{1}(i \in S)$ is the sum of $k$
negatively associated random variables. Therefore, using negative
association, the above expectation is further bounded by
\[
\prod_{i=1}^k
\mathcal{E}\bigl[e^{2na^2\mathbh{1}(i \in S)}\bigr].
\]
Next, for $a$ given by \eqref{EQ:choice_a}, we have
\[
\prod_{i=1}^k\mathcal{E}
\bigl[e^{2na^2\mathbh{1}(i \in S)}\bigr]= \biggl[\bigl(e^{2na^2}-1\bigr)
\frac{k}{p-1}+1 \biggr]^k \le \biggl[1+\frac{1}k
\biggr]^k \le e.
\]

We now show for both cases of the lemma, we have $2a^2k\le1$. Indeed
for case 1, we get $2a^2k= 2(R/4)^{2/q}<2(\log p)/n\le1$. For case 2,
$2a^2k\le1$ follows trivially from the definition of $a$. Also observe
that $k \ge2$ since
\[
\frac{R}{2} \biggl(\frac{\log ({(p-1)}/{4}+1 )}{n} \biggr)^{-q/2}\ge2 \biggl(
\frac{\log p }{n} \biggr)^{q/2} \biggl(\frac{\log
({(p-1)}/{4}+1 )}{n}
\biggr)^{-q/2}\ge2.
\]
This proves part (i) of the statement on $k$. To prove part (ii),
observe that $R^2 < (p-1)n^{-q}/2$ implies that
\[
2^{1-2/q} < 1 < \log \biggl(\frac{p-1}{R^2 n^q}+1 \biggr) ,
\]
which is equivalent to
\[
R n^{q/2} > \frac{R}{2} \biggl(\frac{\log ({(p-1)}/{(R^2 n^q)}+1
)}{2n}
\biggr)^{-q/2}.
\]
Therefore, $k^2 \le R^2 n^q < (p-1)/2$.
\end{pf}

\section{Proof of Theorem~\texorpdfstring{\protect\ref{1595159511599}}{4.1}}
\label{SUB:proof:1595159511599}
The proof follows a similar idea to that of Theorem~\ref{TH:LBquad}.
For the second part of the lower bound, we use exactly the same
construction of two hypotheses as in Lemma~\ref{LEM:chi2div}.
Then it follows that for the $\ell_r$ functional,
\[
\inf_{\hat L} \max_{\Sigma\in\mathcal{H}}
\mathbb{E}_\Sigma \bigl[ \bigl(\hat L-\ell_r(\Sigma)
\bigr)^2 \bigr]\ge Ck^2a^{2r} ,
\]
for some positive constant $C$, where the infimum is taken over all
estimators $\hat L$ of $\ell_r(\Sigma)$ based on $n$ observations. In
case 1, $k^2 a^{2r} = CR^{2r/q}$ while in case 2, following the same
arguments as before,
\[
k^2a^{2r}\ge\frac{R^2}{16} \biggl(\frac{\log (
{(p-1)}/{(R^2n^{q})}+1 )}{2n}
\biggr)^{r-q}\wedge\frac{1}2.
\]
This completes the second part of the lower bound.

The first part of the result is a little bit more complicated than the
construction of $A^{(k)}$ and $B^{(k)}$ in the proof of Theorem~\ref
{TH:LBquad} due to the extra $\log p$ term in the lower bound. We need
to consider a mixture of measures in order to capture the complexity of
the problem. With a slight abuse of notation, we redefine $(2k) \times
(2k)$ matrices $A^{(k)}, B^{(k)}$ as follows:
\[
A^{(k)}=\pmatrix{ \bone
\bone^\top& a \bone\bone^\top
\vspace*{2pt}\cr
a \bone \bone ^\top& \bone\bone^\top } ,\qquad B^{(k)}=\pmatrix{ \bone
\bone^\top& b \bone\bone^\top
\vspace*{2pt}\cr
b \bone \bone ^\top& \bone\bone^\top },
\]
where $a, b \in(0, 1/2)$ and $\bone$ denotes a vector of ones of
length $k$. Since $R^2 < p$, we now construct the block diagonal
covariance matrices
\[
\Sigma_m^{(k)} = \diag(C_1, C_2,
\ldots, C_M, I_{p-2kM}) ,\qquad  m=1,2,\ldots,M,
\]
where the $m$th diagonal block is chosen to be $C_m = B^{(k)}$ while
others are $C_i = A^{(k)}$ for $i \ne m$ and $M = \lfloor p/R \rfloor$.
Also define $\Sigma_0^{(k)}$ to be of the same structure with $C_i =
A^{(k)}$ for all $i$. Then we have $\Sigma_m^{(R/2)} \in\mathcal
{G}_q(R)$ for $m = 0,1,\ldots, M$, since each row of $\Sigma_m^{(R/2)}$
only contains at most $R$ nonzero elements that are bounded by 1.

Let $\mathbb{P}_0$ denote the distribution of $X \sim\mathcal{N}_p(0,
\Sigma_0^{(R/2)})$ and $\mathbb{P}_{m}$ denote the distribution of $X
\sim\mathcal{N}_p(0, \Sigma_m^{(R/2)})$. Let $\mathbb{P}_0^n$ (resp.,
$\mathbb{P}_m^n$) denote the distribution of $\bX=(X_1, \ldots, X_n)$
of $n$ i.i.d. random variables drawn from $\mathbb{P}_0$ (resp.,
$\mathbb{P}_m$). Moreover, let $\bar{\mathbb{P}}^n$ denote the uniform
mixture of $\mathbb{P}_m^n$ over $m \in[M]$. Consider the testing problem
\[
H_0:\qquad \bX\sim\mathbb{P}_0^n\quad \mbox{vs.}\quad
H_1:\qquad \bX\sim\bar{\mathbb{P}}^n.
\]
Using Theorem~2.2, part (iii) of \citet{Tsy09} as before, we need to
show $\chi^2$-divergence can be bounded by a constant. By the same
calculation as in Lemma~\ref{LEM:chi2div}, we have
\[
\chi^2\bigl(\bar{\mathbb{P}}^n, \mathbb{P}_0^n
\bigr)= \frac{1}{M^2} \sum_{1
\le
i,j \le M} \bigl|I-\bigl[\bigl(
\Sigma_0^{(1)}\bigr)^{ -1} \Sigma_i^{(1)}-I
\bigr] \bigl[\bigl(\Sigma _0^{(1)}\bigr)^{ -1}
\Sigma_j^{(1)}-I\bigr]\bigr|^{-n/2} - 1.
\]
Note that $\chi^2$-divergence here depends on $\Sigma_m^{(1)}$ instead
of $\Sigma_m^{(R/2)}$ since perfectly correlated random variables do
not add additional information and hence do not affect $\chi
^2$-divergence (see the proof of Theorem~\ref{TH:LBquad}). Using the
definition of $\Sigma_m^{(1)}$'s, we obtain
\begin{eqnarray*}
& & \bigl|I-\bigl[\bigl(\Sigma_0^{(1)}\bigr)^{ -1}
\Sigma_i^{(1)}-I\bigr] \bigl[\bigl(\Sigma_0^{(1)}
\bigr)^{ -1} \Sigma_j^{(1)}-I\bigr]\bigr|
\\
&&\qquad =  \cases{ %
\displaystyle 1-2\bigl(1+a^2\bigr) \biggl(
\frac{a-b}{1-a^2} \biggr)^2 + \frac{(a-b)^4}{(1-a^2)^2} , & \quad
$\mbox{if } i = j,$
\vspace*{2pt}\cr
1 , &\quad $\mbox{otherwise.}$}
\end{eqnarray*}
Therefore,
\[
\chi^2\bigl(\bar{\mathbb{P}}^n, \mathbb{P}_0^n
\bigr)= \frac{1}{M} \biggl\{ \biggl(1-2\bigl(1+a^2\bigr)
\biggl(\frac{a-b}{1-a^2} \biggr)^2 + \frac
{(a-b)^4}{(1-a^2)^2}
\biggr)^{-n/2} - 1 \biggr\} ,
\]
which is bounded by $( (1-5(a-b)^2)^{-n/2} - 1)/M$ due to the fact
$2(1+a^2)/(1-a^2)^2 \le5$ for $a \le1/2$. Now choose
\[
a = \frac{1}{4} ,\qquad  b = a + \frac{1}{4}\sqrt{\frac{\log p}{n}}.
\]
By assumption, there exists a constant $c_0 >1$ such that $R^2 \le c_0
p$. Thus,
\[
\chi^2\bigl(\bar{\mathbb{P}}^n, \mathbb{P}_0^n
\bigr) \le\frac{1}{M} \biggl\{ \biggl(1-\frac{ \log p}{2 n}
\biggr)^{-n/2} - 1 \biggr\} \le\frac{R}{p} e^{\log
p/2} \le
\sqrt{c_0}.
\]

Using standard techniques to reduce estimation problems to testing
problems as before, we find
\[
\inf_{\hat L} \max_{\Sigma\in\{\Sigma_m^{(R/2)}: m=0,\ldots,M\}}
\mathbb{E}_\Sigma \bigl[ \bigl(\hat L-\ell_r(\Sigma)
\bigr)^2 \bigr]\ge C \bigl(\ell_r\bigl(
\Sigma_0^{(R/2)}\bigr)-\ell_r\bigl(
\Sigma_1^{(R/2)}\bigr) \bigr)^2.
\]
For the above choice of $\Sigma_m^{(k)}$, we have
\[
\bigl(\ell_r\bigl(\Sigma_0^{(R/2)}\bigr)-
\ell_r\bigl(\Sigma_1^{(R/2)}\bigr)
\bigr)^2 = \frac
{R^2}{4} \bigl(b^r -
a^r\bigr)^2 \ge CR^2 \frac{\log p}{n}.
\]
Since $\Sigma_m^{(R/2)} \in\mathcal{G}_q(R)$, the above two displays
imply that
\[
\inf_{\hat L} \max_{\Sigma\in\mathcal{G}_q(R)}
\mathbb{E}_\Sigma \bigl[ \bigl(\hat L-\ell_r(\Sigma)
\bigr)^2 \bigr]\ge C \frac{R^2\log p}{n} ,
\]
which together with the other part of the lower bound, completes the
proof of the theorem.

\section{Proof of Theorem~\texorpdfstring{\protect\ref{35753575357}}{4.4}}
The proof is similar to that of Theorem~\ref{TH:LBquad}, but simpler
since $r \le1$ where no elbow effect exists. Consider hypothesis
construction (\ref{EQ:defSigmaS}) with $\Sigma_S = I_p + \kappa\bar
\Sigma$ and
\[
a = \kappa k^{-1/r} \quad\mbox{and}\quad k = \biggl\lceil R \biggl(
\frac{\log(\nu
p/(R^2 n^q))}{2n} \biggr)^{-q/2} \biggr\rceil.
\]
Choose $\nu$ sufficiently small so that $R(\frac{\log(\nu p/(R^2
n^q))}{2n})^{1-q/2} \le1/2$, which implies $2ka^2 \le1$ and
guarantees the positive semi-definiteness of $\Sigma_S$. Furthermore,
$ka^q \le R$ holds, so $\Sigma_S \in\mathcal G_q(R)$. By the same
derivation as in Theorem~\ref{TH:LBquad}, we are able to show
\[
\chi^2\bigl(\bar{\mathbb{P}}^n, \mathbb{P}_0^n
\bigr) \le e^{\nu} - 1 ,
\]
which by Theorem~2.2(iii) of \citet{Tsy09} leads to the final
conclusion.
\end{appendix}

\begin{supplement}
\stitle{Technical proofs \citet{FanRigWan15}\\}
\slink[doi,text={10.1214/15-\break AOS1357SUPP}]{10.1214/15-AOS1357SUPP} 
\sdatatype{.pdf}
\sfilename{aos1357\_supp.pdf}
\sdescription{This supplementary material contains the introduction to
two-sample high-dimensional testing methods and the proofs of upper
bounds that were omitted from the paper.}
\end{supplement}

%
%
%
%




\printaddresses
\end{document}